\documentclass[preprint]{elsarticle}
\usepackage{graphicx,color}
\usepackage{amsmath,amssymb,amsthm}
\usepackage{subfigure}
\usepackage{lineno,hyperref}
\usepackage{multicol}
\usepackage{color}

\modulolinenumbers[5]

\newtheorem{thm}{Theorem}[section]
\newtheorem{lem}[thm]{Lemma}
\newtheorem{prop}[thm]{Proposition}
\newtheorem{cor}[thm]{Corollary}

\newtheorem{rem}[thm]{Remark}

\begin{document}
%%%%%%%%%%%%%%%%%%%%%%%%%%%%%%%%%%%%%%
\begin{frontmatter}
\title{ Existence of S-shaped type bifurcation curve  with dual cusp catastrophe via  variational methods} 
 
\author{Marcos Leandro Carvalho} 
\ead{marcos$\_$leandro$\_$carvalho@ufg.br}
\address{ Universidade Federal de Goi\'as, Instituto de Matem\'atica\\ 74690-900, Goi\^ania - GO - Brazil}

\author{Yavdat Il'yasov\corref{cor1}}
\ead{ ilyasov02@gmail.com} 
\address{ Institute of Mathematics, Ufa Federal Research Centre, RAS\\
	Chernyshevsky str. 112, 450008 Ufa, Russia\\
	 Universidade Federal de Goi\'as, Instituto de Matem\'atica\\ 74690-900, Goi\^ania - GO - Brazil 
}

\author{Carlos Alberto Santos} 
\ead{csantos@unb.br}
\address{ Universidade de Bras\'ilia, Departamento de Matem\'atica\\
  70910-900, Bras\'ilia - DF - Brazil}

\cortext[cor1]{Corresponding author}

%%%%%%%%%%%%%%%%%%%%%%%%%%%%%%%%%%%%%%
\begin{abstract} We discuss the existence of multiple positive solutions leading to
the occurrence of an S-shaped bifurcation curve to the equations  of the form $
-\Delta_p u= f(\mu,\lambda, u)$
where   $\Delta_p$ is a $p$-Laplacian,  $p>1$,  $\mu, \lambda \in \mathbb{R}$. We deal with relatively unexplored cases when $f(\mu,\lambda, u)$  is non-Lipschitz at $u=0$, $f(\mu,\lambda, 0) = 0 $ and $ f(\mu,\lambda, u) <0$, $u \in (0,r)$, for some $r<+\infty$. We develop the nonlinear generalized Rayleigh quotients method to find a range of parameters where the equation may have distinct branches of positive solutions.
 As a consequence, applying the Nehari manifold method and the mountain pass theorem, we prove that the equation for some range of values $\mu, \lambda$, has at least three positive solutions with two linearly unstable solutions and one  linearly stable.
 The results evidence  that the bifurcation curve is S-shaped and exhibits the so-called dual cusp catastrophe which is characterized by the fact that the corresponding dynamic equation has stable states only within the cusp-shaped region in the control plane of parameters. Our results are new even in the
one-dimensional case and $p=2$.
\end{abstract}

\date{}
\end{frontmatter}

%%%%%%%%%%%%%%%%%%%%%%%%%%%%%%%%%%%%%%%%%%%%%%%%%%%%%%%%%%%%%%%%%%%%%%%%%%%%%%%%%%%%%%%%%%%%%%%%%%%%%%%%%%%%%%55

\section{Introduction}
%Явление гистерезиса возникает вл многих областях и это интересно

Let $\Omega$ be a smooth bounded  domain in $\mathbb{R}^N$, $N\geq 1$. We consider the following  boundary value problem
\begin{equation}\label{p}
\left\{
\begin{aligned}
	-&\Delta_p u= u^{\gamma-1}+\mu u^{\alpha-1}-\lambda u^{q-1}:=f(\mu,\lambda, u)  &&\mbox{in}\ \ \Omega, \\
	&u \geq 0, &&\mbox{in}\ \ \Omega, \\
&u=0                                   &&\mbox{on}\ \ \partial\Omega.
\end{aligned}
\right.
\end{equation}
Here $1<q<\alpha<p<\gamma<p^*$,  $p^*=\frac{pN}{N-p}$ if $p<N$, and $p^*=+\infty$ if $p\geq N$, and $\lambda,\mu \in \mathbb{R}$.

We are interested on the existence of branches of solutions with S-shaped  type bifurcation curve. The S-shaped bifurcation curve arises in the study of many problems such as the Liouville–Bratu–Gelfand equation \cite{gelfand},  the Kolmogorov-Petrovsky-Piscounov equation \cite{Kolmog},  the Minkowski-curvature problem (see, e.g., \cite{Bartnik, Coelho, huang}), in the reaction-diffusion models which described various spatio-temporal phenomena in biology, physics, chemistry, epidemiology, and ecology (see e.g. \cite{Bebernes, Ludw}). This type of bifurcation curve is characterized by  complex and fascinating behavior of systems \cite{Arnold, Gilmore, zeeman}.  The associated dynamic equations could include rarefaction and shock waves (see, e.g.,\cite{Smoller}) and exhibit hysteresis, which is widespread and manifests itself in many physical phenomena, such as dielectric hysteresis, magnetic hysteresis, elastic hysteresis, and some others \cite{Gilmore, Rinzel, Strogatz}.

The "S-shaped bifurcation problem" has
attracted a considerable attention in recent decades  beginning with the celebrated papers by Cohen \cite{cohen} in 1971, Crandall \& Rabinowitz \cite{CrandRabin} in 1973 and Amann \cite{Amann} in 1976. The study of this problem, which usually contains the search for conditions under which the equations have at least three solutions, use different approaches: the bifurcation and continuation methods \cite{CrandRabin, Dancer},  the method of sub-super-solutions \cite{Amann,  Brown, Lee},  the time-map method and the quadrature technique \cite{Brown, huang, Wang}. In \cite{BobDrHer, Lubyshev}, the existence of three positive solution for some a range of parameter of a problem has been obtained by using variational approaches including the mountain pass theorem, the Nehari manifold and fibering methods \cite{Poh}.

Most of the results on the existence of an S-shaped bifurcation curve deal with a nonlinearity $f$ that satisfies $f(u)> 0 $ on $(0, r)$ for some $ r \in (0, +\infty] $, and $f \in C^2[0,+\infty)$ (see e.g. \cite{Amann, Brown,CrandRabin, Lee}). Note that in \eqref{p} we are facing the opposite case, namely, for any $\mu>0$, $\lambda \in \mathbb{R}$, there exists $r_{\mu,\lambda}>0$ such that $f(\mu,\lambda, s)< 0 $,   $s \in (0, r_{\mu,\lambda})$.  Moreover, the nonlinearity in \eqref{p} is  non-Lipschitz at $u=0$ if $p\leq2$ and $\mu\neq 0$ or $\lambda \neq 0$. An additional feature of \eqref{p}, which implies our results  is that the S-shaped bifurcation curve of \eqref{p} exhibits the so-called	dual cusp catastrophe \cite{Gilmore}.

Let us state our main results. The problem \eqref{p} has a variational structure with the energy functional $\Phi_{\lambda,\mu} \in
C^1(W^{1,p}_0)$ given by
\begin{equation*}
\Phi_{\lambda,\mu}(u)=\frac{1}{p}\int |\nabla u|^p+\frac{\lambda}{q}\int |u|^q-\frac{\mu}{\alpha}\int |u|^\alpha-\frac{1}{\gamma}\int  |u|^\gamma.
\end{equation*}
By a weak solution of \eqref{p} we mean a critical point $u$ of $\Phi_{\lambda,\mu}(u)$ on $W^{1,p}_0$, i.e.,
$D\Phi_{\lambda,\mu}(u)=0$. Here and subsequently,  $D\Phi_\lambda(u)$ denotes the Fr\'echet derivative of $\Phi_{\lambda,\mu}$  at $u \in W^{1,p}_0$ and $D\Phi_{\lambda,\mu}(u)(v)$ denotes its the directional derivative  in direction $v  \in W^{1,p}_0$. Hereinafter, for $F \in
C^1(W^{1,p}_0)$ we use the abbreviated notations $F'(tu):=\frac{d}{d t}F(tu)$, $t>0$, $u \in W^{1,p}_0$. 

A nonzero weak
solution $u$ of \eqref{p} is said to be \textit{ground state} if
$\Phi_{\lambda,\mu}(u)\leq \Phi_{\lambda,\mu}(w)$,  
for any  non-zero weak solution $w \in 
W^{1,p}_0$ of \eqref{p}. We say that weak
solution $\bar{u}$ of \eqref{p}  is \textit{linearly stable } if $\bar{u}$ is a local minimizer of $\Phi_{\lambda,\mu}(u)$ in $W^{1,p}_0$ and \textit{linearly unstable} otherwise.

Our approach  is  based on the use of so-called \textit{nonlinear generalized Rayleigh quotient}(\textit{NG-Rayleigh quotient})  whose critical values correspond to the extreme values of the Nehari manifold method \cite{ilyaReil}. 
Following it we introduce the so-called \textit{NG-Rayleigh extremal values} (cf. \cite{MarcCarlIl})
	\begin{align}
&\mu^{e,+}_\lambda:=\inf_{u\in  W^{1,p}_0\setminus 0}\{\mathcal{R}^e_\lambda(u):~(\mathcal{R}_\lambda^e)'(u)=0,~(\mathcal{R}_\lambda^e)''(u)>0\},\label{muplEn}\\
		 &\mu^{e,-}_\lambda:=\inf_{u\in  W^{1,p}_0\setminus 0}\{\mathcal{R}^e_\lambda(u):~(\mathcal{R}_\lambda^e)'(u)=0,~(\mathcal{R}_\lambda^e)''(u)<0\}, \label{muminusEn}\\
		&\mu^{n,+}_\lambda:=\inf_{u\in  W^{1,p}_0\setminus 0}\{\mathcal{R}^n_\lambda(u):~(\mathcal{R}^n_\lambda)'(u)=0,~(\mathcal{R}^n_\lambda)''(u)>0\},\label{muNPlus}\\
					&\mu^{n,-}_\lambda:=\inf_{u\in  W^{1,p}_0\setminus 0}\{\mathcal{R}^n_\lambda(u):~(\mathcal{R}^n_\lambda)'(u)=0,~(\mathcal{R}^n_\lambda)''(u)<0\}. \label{muminus}
	\end{align}
	Here $\mathcal{R}^n_\lambda,\mathcal{R}^e_\lambda:W^{1,p}_0\setminus 0 \to \mathbb{R} $ are the Rayleigh quotients given by
\begin{equation}\label{RaylQS2}
	\begin{aligned}
		&\mathcal{R}^n_\lambda(u):=\frac{\int|\nabla u|^2+\lambda\int |u|^q-\int |u|^\gamma}{\int |u|^\alpha}, ~~u\in W^{1,p}_0\setminus 0,~~ \lambda\in \mathbb{R},\\
		%\label{RaylQS}
		&\mathcal{R}_\lambda^e(u):=\frac{\frac{1}{2}\int|\nabla u|^2+\frac{\lambda}{q}\int |u|^q-\frac{1}{\gamma}\int |u|^\gamma}{\frac{1}{\alpha}\int |u|^\alpha},~~u\in W^{1,p}_0\setminus 0,~~\lambda\in \mathbb{R}. 
	\end{aligned}
\end{equation}
We introduce also the so-called \textit{NG-Rayleigh $\lambda$-extremal value}:   
\begin{align*} 
&\bar{\lambda}=\inf_{u\in  W^{1,p}_0\setminus 0}\frac{(\int|\nabla u|^p)^{\frac{\gamma-q}{\gamma-p}}}{(\int |u|^q) (\int |u|^\gamma)^{\frac{p-q}{\gamma-p}}},\label{lambExtrem}
\end{align*}
and define
\begin{equation}\label{lambdaEN}
	\lambda^e=c^e_{q,\gamma}\bar{\lambda}~~\mbox{and}~ ~\lambda^n=c^n_{q,\gamma}\bar{\lambda},
\end{equation}
where 
$$
	c^e_{q,\gamma}=\frac{q\gamma^{\frac{2-q}{\gamma-2}}}{2^{\frac{\gamma-q}{\gamma-2}}}c^n_{q,\gamma},~~ ~c^n_{q,\gamma}= \frac{(p-\alpha)^{\frac{\gamma-q}{\gamma-p}}(p-q)^{\frac{p-q}{\gamma-q}}(\gamma-p)}{(\alpha-q)(\gamma-\alpha)^{\frac{p-q}{\gamma-p}}(\gamma-q)^{\frac{\gamma-q}{\gamma-p}}}.
$$

\begin{lem} \label{lem1}
Assume that $1<q<\alpha<p<\gamma<p^*$. Then
	\begin{description}
		\item[($1^o$)] 
		 $ 0<\lambda^e<\lambda^n<+\infty$,
		\item[($2^o$)]   $0<\mu^{n,+}_\lambda<   \mu^{e,+}_\lambda<\mu^{e,-}_\lambda<\mu^{n,-}_\lambda<+\infty$, for  $\lambda \in  (0,\lambda^e)$.
		\end{description}
\end{lem}
%$(u^0_{\lambda, \mu}), (u^1_{\lambda, \mu}), (u^2_{\lambda, \mu})$
%Let us state our main results. 

We deal with the Nehari manifold 
\begin{equation*}\label{eq:NehariR}
		\mathcal{N}_{\lambda,\mu}=\{u \in W^{1,p}_0\setminus 0: ~\Phi_{\lambda,\mu}'(u)=0\}.
\end{equation*}
A local minimum or maximum point of the function $\Phi_{\lambda,\mu}(u)$ subject to $\mathcal{N}_{\lambda,\mu}$ is called the extremal point of $\Phi_{\lambda,\mu}(u)$ on the Nehari manifold. We denote
$$
\hat{\Phi}_{\lambda,\mu}=\min_{u \in \mathcal{N}_{\lambda, \mu}}\Phi_{\lambda,\mu}(u)~~\mbox{and}~~\mathcal{M}_{\lambda, \mu}:=\{u \in \mathcal{N}_{\lambda,\mu}:~~ \hat{\Phi}_{\lambda,\mu}=\Phi_{\lambda,\mu}(u) \}.
$$

\begin{thm}\label{thm1}
	Assume that $1<q<\alpha<p<\gamma<p^*$,  $\lambda \in (0,\lambda^e)$.  Then  problem \eqref{p} admits two distinct branches of positive solutions  $(u^2_{\lambda,\mu})$, $(u^3_{\lambda,\mu})$, namely:   
\begin{description}
\item[{\rm $(1^0)$}]  There exists $\hat{\mu}_\lambda^{*} \in (\mu^{n,+}_\lambda, \mu^{e,+}_\lambda)$ such that for any $\mu\in [\hat{\mu}_\lambda^{*},\mu^{n,-}_\lambda)$, problem \eqref{p} possesses  a positive  solution $u^2_{\lambda, \mu}\in C^{1,\kappa}(\overline{\Omega})$, $\kappa \in (0,1)$ 
   such that (i) $u^2_{\lambda, \mu}$ is a ground state; (ii) $u^2_{\lambda,\mu}$ is linearly stable and  $\Phi''_{\lambda,\mu 	}(u^2_{\lambda, \mu})>0$; (iii) the function  $ \mu \mapsto \Phi_{\lambda,\mu}(u^2_{\lambda, \mu})$ is continuous and monotone decreasing on $(\hat{\mu}_\lambda^{*},\,\mu^{n,-}_\lambda)$; (iv)  $\Phi_{\lambda,\mu}(u^2_{\lambda, \mu})>0$ for $\mu\in (\hat{\mu}_\lambda^{*},\mu^{e,+}_\lambda)$, $\Phi_{\lambda,\mu^{e,+}_\lambda}(u^2_{\lambda, \mu^{e,+}_\lambda})=0$, and $\Phi_{\lambda,\mu}(u^2_{\lambda, \mu})<0$ for $\mu\in (\mu^{e,+}_\lambda,\mu^{n,-}_\lambda)$. 
	
\item[{\rm $(2^0)$}] For any  $\mu \in (-\infty, \mu^{n,-}_\lambda)$, problem \eqref{p} possesses  a positive solution $u^3_{\lambda,\mu}\in C^{1,\kappa}(\overline{\Omega})$, $\kappa \in (0,1)$ such that (i) $u^3_{\lambda,\mu}$ is linearly unstable and $\Phi''_{\lambda,\mu
	}(u^3_{\lambda, \mu})<0$; (ii) the function  $(-\infty, \mu^{n,-}_\lambda) \ni \mu \mapsto \Phi_{\lambda,\mu}(u^3_{\lambda, \mu})$ is continuous and monotone decreasing;  (iii)  $\Phi_{\lambda,\mu}(u^3_{\lambda, \mu})>0$ if $\mu \in (-\infty, \mu^{e,-}_\lambda)$, $\Phi_{\lambda,\mu^{e,-}_\lambda}(u^3_{\lambda, \mu^{e,-}_\lambda})=0$,  and  $\Phi_{\lambda,\mu}(u^3_{\lambda, \mu})<0$ if $\mu \in (\mu^{e,-}_\lambda,\mu^{n,-}_\lambda)$; (iv)  $\Phi_{\lambda,\mu}(u^3_{\lambda, \mu}) \to +\infty$, $\|u^3_{\lambda, \mu}\|_1 \to +\infty$ as $\mu \to -\infty$. 
	
	Moreover, if $\mu \in (-\infty, \mu^{n,+}_\lambda)$, then $u^3_{\lambda,\mu}$ is a ground state of \eqref{p}.
\end{description}

\end{thm}
A weak solution $u_{\lambda,\mu} \in W^{1,p}_0\setminus 0$ of \eqref{p} is said to be mountain pass type if  
\begin{equation}\label{MOUNT}
	\Phi_{\lambda,\mu}(u_{\lambda,\mu})=\hat{\Phi}_{\lambda,\mu}^{m}:=\inf_{\gamma\in P}\sup_{u\in \gamma} \Phi_{\lambda,\mu}(u)>0,
\end{equation} 
for the paths set $P:=\{g \in C([0,1];W^{1,p}_0):g(0)=0,~ g(1)=w_1\}$ with some $w_1 \in W^{1,p}_0$ such that $\Phi_{\lambda,\mu}(w_1)<0$.

\begin{thm}\label{thm2}
	Assume that $1<q<\alpha<p<\gamma<p^*$, $\lambda>0$, $-\infty< \mu <+\infty$. Then \eqref{p} admits a positive  mountain pass type solution $u_{\lambda,\mu}\in C^{1,\kappa}(\overline{\Omega})$, $\kappa \in (0,1)$ such that 	$\Phi_{\lambda,\mu}(u_{\lambda,\mu})>0$. Moreover, if $\lambda \in (0,\lambda^e)$, $\mu \in (-\infty, \mu^{n,+}_\lambda)$, then (i)  $u_{\lambda,\mu}$ is a ground state of \eqref{p}, i.e., $u_{\lambda,\mu} \in \mathcal{M}_{\lambda, \mu}$; (ii)  $\Phi_{\lambda,\mu}(u_{\lambda, \mu})\to +\infty$, $\|u_{\lambda, \mu}\|_1 \to +\infty~~\mbox{as}~~\mu \to -\infty$, (iii)  $\Phi_{\lambda,\mu}(u_{\lambda, \mu})\to 0$, $\|u_{\lambda, \mu}\|_1 \to 0$ as $\mu \to +\infty$.
 \end{thm}

\begin{figure}[ht]
	\centering
	\includegraphics[width=0.5\linewidth]{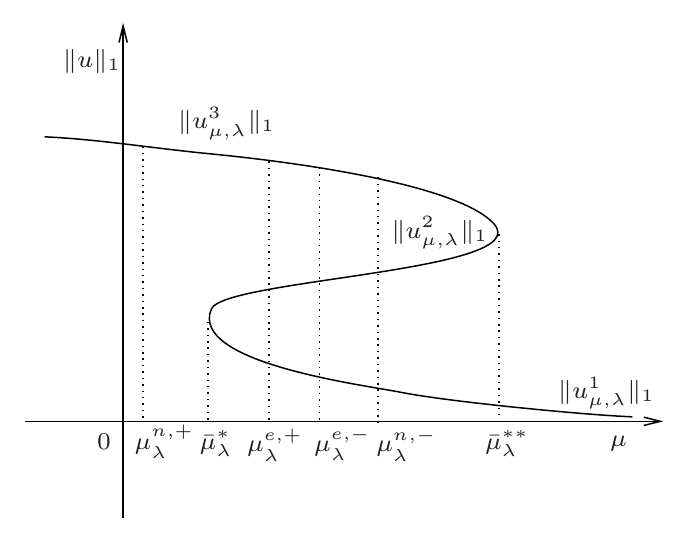}\\
	\caption{{\it The branches  of solutions to \eqref{p} in term of  the norm $\|\cdot\|_1$}}
	\label{fig4}
\end{figure}

\begin{figure}[ht]
	\centering
	\includegraphics[width=0.5\linewidth]{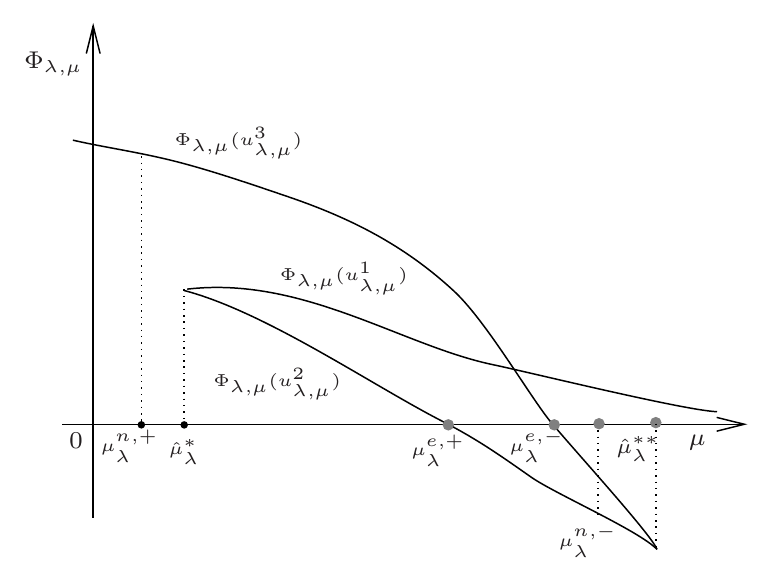}\\
	\caption{{\it The branches of   solutions to \eqref{p} in term of the energy levels $\Phi_{\lambda,\mu}$.}}
	\label{fig5}
\end{figure}

Theorems  \ref{thm1}, \ref{thm2} yield the following result on the existence of three distinct branches of weak positive solutions of \eqref{p}. 
\begin{thm}\label{thm3}
Assume $\lambda \in (0,\lambda^e)$ and $\mu\in[\mu^{e,-}_\lambda,\mu_\lambda^{n,-})$. Then \eqref{p} admits at least three distinct positive solutions: $u_{\lambda,\mu}^1$, $u^2_{\lambda,\mu}$, $u^3_{\lambda,\mu}$ such that $\Phi_{\lambda,\mu}(u_{\lambda,\mu}^1)>0$, $\Phi_{\lambda,\mu}(u^2_{\lambda, \mu})<0$, $\Phi_{\lambda,\mu}(u^3_{\lambda, \mu})\leq 0$ and $\Phi''_{\lambda,\mu}(u^2_{\lambda, \mu})>0$, $\Phi''_{\lambda,\mu}(u^3_{\lambda, \mu})<0$. Furthermore, $u^2_{\lambda,\mu}$ is linearly stable while  $u^1_{\lambda,\mu}, u^3_{\lambda,\mu}$ are   linearly unstable solutions.
\end{thm}
In view of these result, it is natural to expect that the branches of solutions to \eqref{p} behave as depicted
on Figures \ref{fig4}, \ref{fig5}. That is  there are two bifurcation values  $\bar{\mu}_\lambda^{*}$, 
$\bar{\mu}_\lambda^{**}$ such that: (i) $0<\mu^{n,+}_\lambda\leq \bar{\mu}_\lambda^{*}\leq \hat{\mu}^*_\lambda \leq \mu^{e,+}_\lambda<\mu^{n,-}_\lambda\leq \bar{\mu}_\lambda^{**}$; (ii) the branches  $u^1_{\lambda, \mu}$, $u^2_{\lambda, \mu}$ continue until the bifurcation value $\hat{\mu}^{*}_\lambda$, where they coincide; (iii)  $ u^2_{\lambda, \mu} $, $ u^3_{\lambda, \mu} $ continue until the bifurcation value $ \hat{\mu}_\lambda^{**}$, where they coincide. Furthermore, we anticipate that the solution $u^2_{\lambda, \mu}$ is stable for any $\mu \in (\hat{\mu}^{*}_\lambda, \bar{\mu}_\lambda^{**})$, whereas 
solutions $u^1_{\lambda, \mu}$ for $\mu \in (-\infty, \bar{\mu}^{*}_\lambda)$ and $u^3_{\lambda, \mu}$ for $\mu \in (\bar{\mu}_\lambda^{**}, +\infty)$ are unstable. It is important to emphasize that such behavior means that the S-shaped bifurcation curve of \eqref{p} exhibits the so-called dual cusp catastrophe \cite{Gilmore}. This type of catastrophe is characterized by the fact that the corresponding dynamic equation has an opposite behaviour, namely it has stable states only within the cusp-shaped region in the control plane of parameters $\mu, \lambda$. Note that the cusp catastrophe, which is more common in the studying S-shaped bifurcation curves (see, e.g., \cite{BobDrHer, Brown,  Korman, Wang}) has a stable state for the entire range of parameters and  are characterized by hysteresis behaviour \cite{Gilmore, Rinzel, Strogatz}.  We failed to find any reference addressing the existence of an S-shaped bifurcation curve of nonlinear PDE's in high dimensions which exhibits the dual cusp catastrophe.

%In particular, we could not find   any reference addressing the existence of S-shaped bifurcation curve  with dual cusp catastrophe  to PDEs.

%
%Notice that the usual cusp catastrophe encountered in the study of S-shaped bifurcation curves is globally stable, that is, it  has a stable state for an entire range of parameters, and the system is characterized by hysteresis behaviour \cite{Gilmore, Rinzel, Strogatz}.

\begin{rem}
	We anticipate that $\bar{\mu}_\lambda^{*}=\hat{\mu}^*_\lambda$, and $u^3_{\lambda, \mu}$ is a ground state for $\mu \in (-\infty, \bar{\mu}_\lambda^{*})$; $u^2_{\lambda, \mu}$ is a ground state for $\mu \in 
		(\bar{\mu}^{*}_\lambda, \bar{\mu}_\lambda^{**})$; $u^1_{\lambda, \mu}$ is a ground state for $\mu \in (\bar{\mu}_\lambda^{**}, -\infty)$. In this regard, it is interesting to note such a phenomenon for the ground states branch as the appearance of jump at $\bar{\mu}^{*}_\lambda$ of the energy level $\Phi_{\lambda,\mu}$  and jump at  $\bar{\mu}_\lambda^{**}$ of the values of norm $\|\cdot\|_1$ (see Figures \ref{fig4}, \ref{fig5})
\end{rem}
.

\begin{rem}
    Evidently, if $u$ is a solution of \eqref{p}, then so is $-u$.
    Thus, Theorems \ref{thm1}-\ref{thm3}
    actually establish the existence at least of pairs of the branches of solutions $(u^1_{\lambda,\mu},-u^1_{\lambda,\mu})$, $(u^2_{\lambda,\mu},-u^2_{\lambda,\mu})$, and $(u^3_{\lambda,\mu},-u^3_{\lambda,\mu})$, respectively.
\end{rem}

\begin{rem}
The present article extends  results obtained in our previous work \cite{MarcCarlIl}, where   \eqref{p}  has been studied  in the case $p = 2$, and the solution $u^2_{\lambda, \mu}$ had been obtained for $\lambda \in (0,\lambda^e)$, $\mu \in (\mu^{e,+}_\lambda,\mu_\lambda^{n,-})$,  whereas the solution $u^3_{\lambda, \mu}$ had been obtained, under  additional restrictions $1+\alpha<\gamma<2^*$, for $\lambda \in (0,\lambda^n)$, $\mu \in (-\infty,\mu_\lambda^{n,-})$. Furthermore, we succeed to develop a new approach, which proved to be useful in obtaining solutions on wider intervals of parameters and improving the results in general.
\end{rem}

\section{ Preliminaries}\label{sec2}

Hereinafter, we use the standard notation $L^r:=L^r(\Omega
)$, $1 \leq r\leq +\infty$ for the Lebesgue spaces endowed with the norm  $\|u\|_{L^r}:=(\int |u|^r)^{1/r}$,  $W^{1,p}_0:=W^{1,p}_0(\Omega)$ for the Sobolev space endowed with the norm  $\|u\|_1:=(\int |\nabla u|^p)^{1/p}$. The weak convergence in $W_{0}^{1,p}$ we shall denote by "$\rightharpoonup $".

\begin{lem}\label{lemNL}
 Let $\hat{u} \in \mathcal{N}_{\lambda,\mu}$ be an extremal point  of $\Phi_{\lambda,\mu}(u)$ on the Nehari manifold.  Suppose that $\Phi_{\lambda,\mu}''(\hat{u}):=D\Phi_{\lambda,\mu}'(\hat{u})(\hat{u}) \neq 0$. Then $D\Phi_{\lambda,\mu}(\hat{u})=0$.
	\end{lem}
\begin{proof} Due to the assumption we may apply the Lagrange multiplier rule   (see  Proposition 43.19 in \cite{zeidl}), and thus, we have $D\Phi_{\lambda,\mu}(\hat{u})+\nu D\Phi_{\lambda,\mu}'(\hat{u})=0$, for some $\nu \in \mathbb{R}$. 
Testing this equality by $\hat{u}$ we obtain $\mu D\Phi_{\lambda,\mu}'(\hat{u})(\hat{u})=0$. Since $D\Phi_{\lambda,\mu}'(\hat{u})(\hat{u}) \neq 0$, $\nu=0$, and therefore, $D\Phi_{\lambda,\mu}(\hat{u})=0$.
\end{proof}
Observe that 	$\Phi_{\lambda,\mu}$ is coercive on 
		$\mathcal{N}_{\lambda,\mu}$, $\forall \lambda >0$, $\forall\mu \in \mathbb{R}$. Indeed, by the Sobolev inequality,
\begin{align*}\label{Phi-bounded}
							\Phi_{\lambda,\mu}(u)\geq \frac{\gamma-p}{p\gamma}\|u\|_1^p-\mu C\|u\|_1^\alpha,~~\forall u\in \mathcal{N}_{\lambda,\mu},
							\end{align*}
for some  constant $C>0$  independent of $u \in \mathcal{N}_{\lambda,\mu}$. Since $\alpha<p$, we have   $\Phi_{\lambda,\mu}(u) \to +\infty$ for  $\|u\|_1 \to +\infty$ and $u \in \mathcal{N}_{\lambda,\mu}$.

\begin{figure}[!ht]
\begin{minipage}[h]{0.49\linewidth}
\center{\includegraphics[scale=0.7]{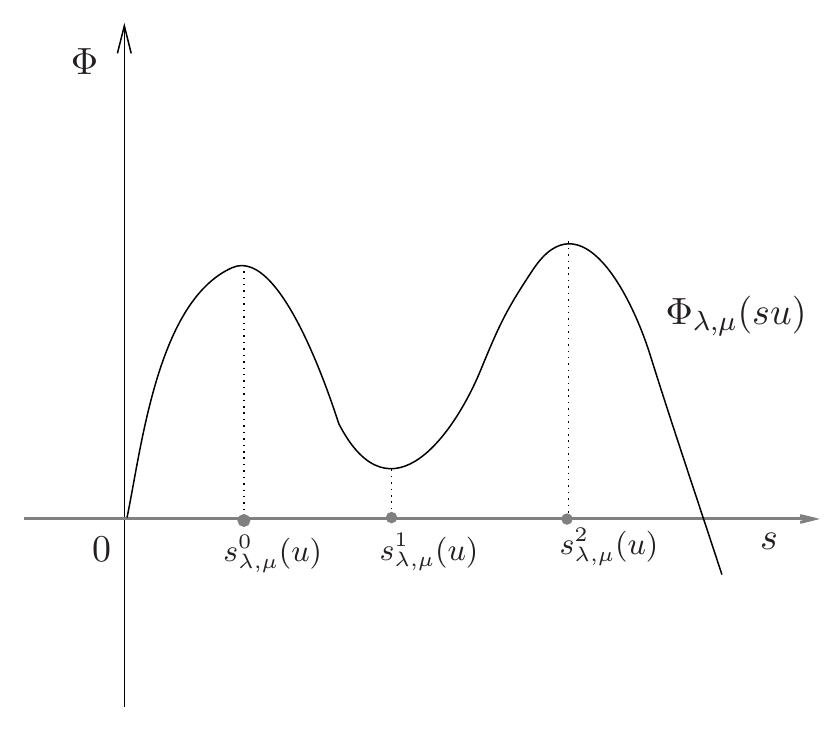}}
\caption{{\it Fibering function $\Phi_{\lambda,\mu}(su)$, $s\geq 0$,  $u \in W^{1,p}_0 $.}}
	\label{fig1}
\end{minipage}
\hfill
\begin{minipage}[h]{0.49\linewidth}
\center{\includegraphics[scale=0.7]{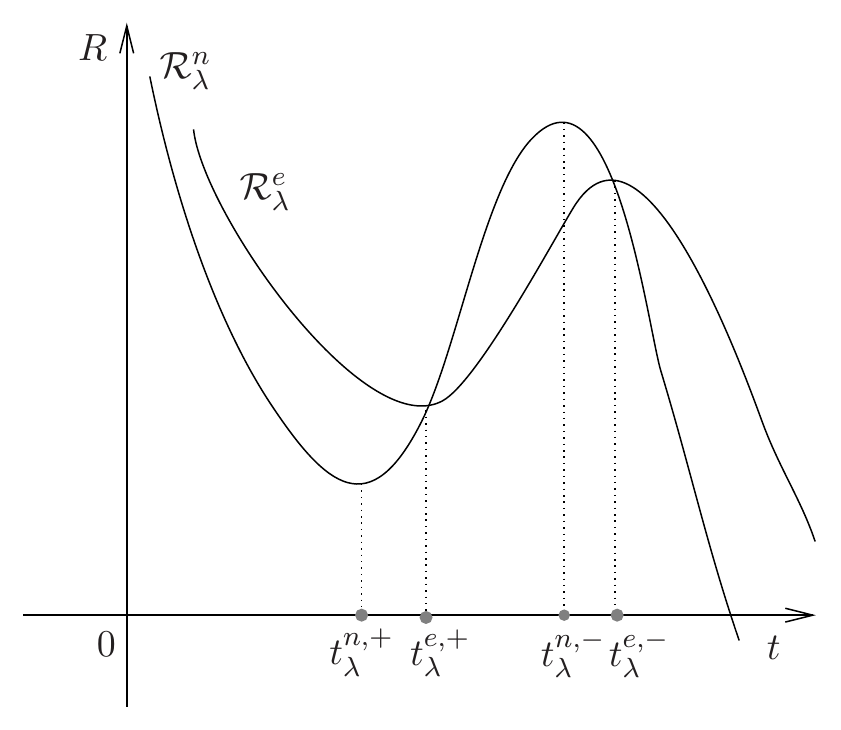}}
\caption{\textit{The functions $\mathcal{R}^e_\lambda(tu)$, $\mathcal{R}^n_\lambda(tu)$}}
\label{fig2}
\end{minipage}
\end{figure}

It is easily seen that for $u\in  W^{1,p}_0\setminus 0$, the fibering function $\Phi_{\lambda,\mu}(su)$, $s>0$,  may have at most three nonzero critical points 
$$0< s^1_{\lambda,\mu}(u)\leq  s^2_{\lambda,\mu}(u) \leq s^3_{\lambda,\mu}(u)<\infty
$$ 
such that 
$
\Phi''_{\lambda,\mu}(s^1_{\lambda,\mu}(u) u)\leq 0,~~~ \Phi''_{\lambda,\mu}(s^2_{\lambda,\mu}(u) u)\geq 0 ,~~~ \Phi''_{\lambda,\mu}(s^3_{\lambda,\mu}(u) u)\leq 0
$ (see Figure \ref{fig1}).

To apply the Nehari manifold method, we need to find values $\lambda$, $\mu$, where the strong inequalities $0< s^1_{\lambda,\mu}(u)< s^2_{\lambda,\mu}(u) <s^3_{\lambda,\mu}(u)$ hold.   
We solve this  by the recursively application of the nonlinear generalized Rayleigh quotient method  proposed in \cite{MarcCarlIl}. In the first step of this recursive procedure, we consider
\begin{equation}\label{RaylQS}
	\mathcal{R}^n_\lambda(u)=\frac{\int|\nabla u|^p+\lambda\int |u|^q-\int |u|^\gamma}{\int |u|^\alpha}, ~~u\in  W^{1,p}_0\setminus 0, ~\lambda\in \mathbb{R}. 
\end{equation}
Notice that for $u\in  W^{1,p}_0\setminus 0$,
\begin{equation}\label{RandPhi}
\begin{aligned}
	&\mathcal{R}^n_\lambda(u)=\mu~	\Leftrightarrow~\Phi'_{\lambda,\mu}(u)=0\\
	&\mathcal{R}^n_\lambda(u)=\mu,~ (\mathcal{R}^n_\lambda)'(u)>(<)\, 0~	\Leftrightarrow~\Phi''_{\lambda,\mu}(u)>(<)\, 0.
\end{aligned}
\end{equation}
In particular,
$$
\mathcal{N}_{\lambda,\mu}=\{u \in W^{1,p}_0\setminus 0: ~\mathcal{R}^n_\lambda(u)=\mu\}.
$$
Moreover, since for any $\lambda>0$,  $\mathcal{R}^n_\lambda(tu) \to +\infty$ as $t\downarrow 0$ and 
$\mathcal{R}^n_\lambda(tu) \to -\infty$ as $t\to +\infty$,
\begin{equation}\label{neqemptysetN}
	\mathcal{N}_{\lambda,\mu}\neq \emptyset,~~~\forall \mu \in (-\infty,+\infty),~ \forall \lambda\in (0,+\infty).
\end{equation}

Using the Implicit Function Theorem we have
\begin{prop}\label{ImplS} 
If  $(\mathcal{R}^n_\lambda)'(s^i_{\lambda,\mu}(u_0)u_0)\neq 0$ for  $u_0 \in W^{1,p}_0\setminus 0$, $i=1,2,3$, then  there exists a neighbourhood  $U_{u_0}\subset W^{1,p}_0\setminus 0$ of $u_0$  such that $s^i_{\lambda,\mu}(\cdot) \in C^1(U_{u_0})$.
\end{prop}

Simple analysis shows that for any given $u \in W^{1,p}_0\setminus 0$, the fibering function $ \mathcal{R}^n_\lambda(tu)$ may have at most two non-zero  critical points $t_\lambda^{n,+}(u), t_\lambda^{n,-}(u)$,  where   $t_\lambda^{n,+}(u)$ is a local minimum, $t_\lambda^{n,-}(u)$ is a local maximum point of $ \mathcal{R}^n_\lambda(tu)$, and 
\begin{equation}\label{tsts}
	0<s^1_{\lambda,\mu}(u)\leq t_\lambda^{n,+}(u) \leq s^2_{\lambda,\mu}(u)\leq t_\lambda^{n,+}(u) \leq s^3_{\lambda,\mu}(u)<\infty.
\end{equation}
(see Figure \ref{fig2}).

To split the points $t_\lambda^{n,+}(u), t_\lambda^{n,-}(u)$,  in the second step of the recursive procedure, we apply the nonlinear generalized Rayleigh quotient method to the functional $\mathcal{R}^n_\lambda$  with respect to the parameter  $\lambda$, i.e., we consider 
\begin{equation}\label{110}
\Lambda^n(u):= \frac{(p-\alpha) \int|\nabla u|^p-(\gamma-\alpha)\int |u|^\gamma}{(\alpha-q)\int |u|^q}, ~~u \in W^{1,p}_0\setminus 0.
\end{equation}
Notice that for any $u \in W^{1,p}_0\setminus 0$, $(\mathcal{R}^n_\lambda)'(tu)=0$
$\Leftrightarrow$ $\Lambda^n(tu)=\lambda$.
The only solution of
$
\displaystyle{\frac{d}{dt}\Lambda^n(tu)=0}
$
is a   global maximum point $t^n(u)$ of the function $\Lambda^n(tu)$ which can be found precisely   
\begin{equation}
\label{17}
t^n(u):= \left(C_n\frac{\int|\nabla u|^p}{\int |u|^\gamma}\right)^{1/(\gamma-p)},~~\forall u\in  W^{1,p}_0\setminus 0,
\end{equation}
where 
$$
C_n=\frac{(p-\alpha)(p-q)}{(\gamma-\alpha)(\gamma-q)}.
$$ 
This allow us to introduce the following  \textit{NG-Rayleigh  $\lambda$-quotient}   
\begin{equation}\label{const0}
\lambda^n(u):=\Lambda^n(t^n(u)u)=c^n_{q,\gamma}\frac{(\int|\nabla u|^p)^{\frac{\gamma-q}{\gamma-p}}}{(\int |u|^q) (\int |u|^\gamma)^{\frac{p-q}{\gamma-p}}},
\end{equation}
and the corresponding  \textit{ principal extremal value}   
\begin{align*}\label{lambExtrem0}
&\lambda^n=\inf_{u\in   W^{1,p}_0\setminus 0}\sup_{t>0}\Lambda^n(tu),
=	 c^n_{q,\gamma}\frac{(\int|\nabla u|^p)^{\frac{\gamma-q}{\gamma-p}}}{(\int |u|^q) (\int |u|^\gamma)^{\frac{p-q}{\gamma-p}}},
\end{align*}
where 
$$
	c^n_{q,\gamma}= \frac{(p-\alpha)^{\frac{\gamma-q}{\gamma-p}}(p-q)^{\frac{p-q}{\gamma-q}}(\gamma-p)}{(\alpha-q)(\gamma-\alpha)^{\frac{p-q}{\gamma-p}}(\gamma-q)^{\frac{\gamma-q}{\gamma-p}}}.
$$
Note that this definition of $\lambda^n$ coincides with  \eqref{lambdaEN}.

It easily follows  (cf. \cite{MarcCarlIl})
\begin{prop}\label{propRN2} 
	For any $\lambda\in (0,\lambda^n)$ and $u\in W^{1,p}_0\setminus 0$, the function $\mathcal{R}^n_\lambda(tu)$ has precisely two distinct critical points
such that $0<t_\lambda^{n,+}(u)<t_\lambda^{n,-}(u)$,  with $t_\lambda^{n,+}(\cdot), t_\lambda^{n,-}(\cdot) \in C^1(W^{1,p}_0\setminus 0)$. Moreover,  
\begin{itemize}
	\item $(\mathcal{R}^n_\lambda)''(t_\lambda^{n,+}(u)u)>0$,  $(\mathcal{R}^n_\lambda)''(t_\lambda^{n,-}(u)u)<0$,
	\item $(\mathcal{R}^n_\lambda)'(tu)<0$ $\Leftrightarrow$ $t \in (0,t_\lambda^{n,+}(u)) \cup (t_\lambda^{n,-}(u),\infty)$,
	\item $(\mathcal{R}^n_\lambda)'(tu)>0$ $\Leftrightarrow$ $t \in (t_\lambda^{n,+}(u),t_\lambda^{n,-}(u))$.
\end{itemize}
 \end{prop}
Observe that this and \eqref{tsts} imply that $0<s^1_{\lambda,\mu}(u)<s^3_{\lambda,\mu}(u)<\infty$ for any $\lambda\in (0,\lambda^n)$, $u\in   W^{1,p}_0\setminus 0$. Thus, for $\lambda\in (0,\lambda^n)$, we are able to introduce the following \textit{NG-Rayleigh $\mu$-quotients}
\begin{equation*}\label{mumuN}
\mu^{n,+}_\lambda(u):=\mathcal{R}^n_\lambda(t_\lambda^{n,+}(u)u),~~ \mu^{n,-}_\lambda(u):=\mathcal{R}^n_\lambda(t_\lambda^{n,-}(u)u), ~~u\in  W^{1,p}_0\setminus 0.
\end{equation*}
By Proposition \ref{propRN2} and regularity of $\mathcal{R}^n_\lambda$ it follows that $\mu^{n,+}_\lambda$ and $\mu^{n,-}_\lambda$ are  $C^1(W^{1,p}_0\setminus 0)$,  $\lambda\in (0,\lambda^n)$.
It is easily seen that the corresponding \textit{principal extremal values} 
	\begin{align}
		&\mu^{n,+}_\lambda=\inf_{u\in   W^{1,p}_0\setminus 0}\mu^{n,+}_\lambda(u), \label{mupl}\\
		&\mu^{n,-}_\lambda=\inf_{u\in   W^{1,p}_0\setminus 0}\mu^{n,-}_\lambda(u)\label{muminusSH}
	\end{align}
coincide with \eqref{muNPlus} and \eqref{muminus}, respectively.

 We also need  the so-called zero-energy level   Rayleigh quotient
\begin{align*}
	&\mathcal{R}_\lambda^e(u)=\frac{\frac{1}{p}\int|\nabla u|^p+\frac{\lambda}{q}\int |u|^q-\frac{1}{\gamma}\int |u|^\gamma}{\frac{1}{\alpha}\int |u|^\alpha},~~u\in W^{1,p}_0\setminus 0, 
	\end{align*}
	which is characterized by the fact that $\mathcal{R}^e_\lambda(u)=\mu$ $\Leftrightarrow$ $\Phi_{\lambda,\mu}(u)=0$. It is easy to see that  $\mathcal{R}_\lambda^e(u)$ possesses similar properties to that $\mathcal{R}_\lambda^n(u)$. In particular, the fibering function $ \mathcal{R}^e_\lambda(tu)$ may have at most two non-zero fibering critical points $0<t_\lambda^{e,+}(u)\leq t_\lambda^{e,-}(u)<+\infty$ so that $t_\lambda^{e,+}(u)$ is a local minimum  while $t_\lambda^{e,-}(u)$ is a local maximum point of $ \mathcal{R}^e_\lambda(tu)$. Moreover, the same conclusion as for $\Lambda^n(u)$ can be drawn for the Rayleigh quotient 
\begin{equation}\label{11E}
\Lambda^e(u):= q\frac{\frac{(p-\alpha)}{p} \int|\nabla u|^p-\frac{(\gamma-\alpha)}{\gamma}\int |u|^\gamma}{(\alpha-q)\int |u|^q},
\end{equation}
which is characterized by the fact that   $(\mathcal{R}^e_\lambda)'(tu)=0$ $\Leftrightarrow$ $\Lambda^e(tu)=0$ for any $u \in W^{1,p}_0\setminus 0$.
The unique solution of
$
\frac{d}{dt}\Lambda^e(tu)=0
$
is a  global maximum point of the function $\Lambda^e(tu)$ defined by  
\begin{equation}
\label{177}
t^e(u):= \left(C_e\frac{\| u\|^p_1}{\| u\|^\gamma_{L^\gamma}}\right)^{1/(\gamma-p)},~~\forall u\in  W^{1,p}_0\setminus 0,
\end{equation}
where 
$$
C_e=\frac{\gamma(p-\alpha)(p-q)}{p(\gamma-\alpha)(\gamma-q)}.
$$ 
Thus we have the following  NG-Rayleigh  quotient $\lambda^e(u):= \Lambda^e(t^e(u)u)$, $u\in  W^{1,p}_0\setminus 0$ with the corresponding principal extremal value   
\begin{align*} 
			 &\lambda^e=\inf_{u\in  W^{1,p}_0\setminus 0}\sup_{t>0}\Lambda^e(tu)=c^e_{q,\gamma}\inf_{u\in  W^{1,p}_0\setminus 0}\frac{\| u\|_1^{p\frac{\gamma-q}{\gamma-p}}}{\| u\|^q_{L^q} \| u\|_{L^\gamma}^{\gamma\frac{p-q}{\gamma-p}}}. \label{lambExtrem1}
	\end{align*}
Note that this definition of $\lambda^e$ coincides with  \eqref{lambdaEN}.
We thus have
\begin{prop}\label{propRNE2} 
	For any $\lambda\in (0,\lambda^e)$ and $u\in W^{1,p}_0\setminus 0$, the function $\mathcal{R}^e_\lambda(tu)$ has precisely two distinct critical points such that $0<t_\lambda^{e,+}(u)<t_\lambda^{e,-}(u)$ with $t_\lambda^{e,+}(\cdot), t_\lambda^{e,-}(\cdot) \in C^1(W^{1,p}_0\setminus 0)$. Moreover,  
\begin{itemize}
	\item $(\mathcal{R}^e_\lambda)''(t_\lambda^{e,+}(u)u)>0$ and $(\mathcal{R}^e_\lambda)''(t_\lambda^{e,-}(u)u)<0$,
	\item $(\mathcal{R}^e_\lambda)'(tu)<0$ $\Leftrightarrow$ $t \in (0,t_\lambda^{e,+}(u)) \cup (t_\lambda^{e,-}(u),\infty)$,
	\item $(\mathcal{R}^e_\lambda)'(tu)>0$ $\Leftrightarrow$ $t \in (t_\lambda^{e,+}(u),t_\lambda^{e,-}(u))$.
\end{itemize}
\end{prop} 
Hence, for $\lambda\in (0,\lambda^e)$, we are able to introduce the following zero-energy level \textit{NG-Rayleigh $\mu$-quotients}
\begin{equation*}\label{mumuN4}
\mu^{e,+}_\lambda(u):=\mathcal{R}^e_\lambda(t_\lambda^{e,+}(u)u),~~~~ \mu^{e,-}_\lambda(u):=\mathcal{R}^n_\lambda(t_\lambda^{e,-}(u)u), ~~u\in  W^{1,p}_0\setminus 0.
\end{equation*}
It is easily seen that the corresponding zero-energy level  principal extremal values
	\begin{align}\label{mumue}
		&\mu^{e,+}_\lambda=\inf_{u\in   W^{1,p}_0\setminus 0}\mu^{e,+}_\lambda(u),~~\mu^{e,-}_\lambda=\inf_{u\in   W^{1,p}_0\setminus 0}\mu^{e,-}_\lambda(u)
	\end{align}
	coincide with \eqref{muplEn} and \eqref{muminusEn}, respectively.

The relationships among the above introduced Rayleigh quotients are given by the following lemma (see Figures \ref{fig2})

\begin{lem}\label{lemGEOM}
Assume that $1<q<\alpha<p<\gamma$, $u\in W^{1,p}_0\setminus 0$, $t>0$. 
		\begin{enumerate}
			%\item[\rm{$(i)$}] $\Lambda^e(tu)<\Lambda^n(tu)$ for   sufficiently small $t>0$,
			\item[\rm{$(i)$}] $\Lambda^{e}(tu)=\Lambda^n(tu)$   $\Leftrightarrow $ $t=t^e(u)$,
			\item[\rm{$(ii)$}]    $\mathcal{R}_\lambda^e(tu)=\mathcal{R}_\lambda^n(tu)$ $\Leftrightarrow $ $t=t_\lambda^{e,+}(u)$ or $t=t_\lambda^{e,-}(u)$, $\forall \lambda \in(0,\lambda^e)$,
			\item[\rm{$(iii)$}]   $t_\lambda^{n,+}(u)<t_\lambda^{e,+}(u)<t^e(u)<t_\lambda^{n,-}(u)<t_\lambda^{e,-}(u)$, $\forall \lambda \in(0,\lambda^e)$,
			\item[\rm{$(iv)$}]  $\mathcal{R}_\lambda^n(tu)<\mathcal{R}_\lambda^e(tu)$ $\Leftrightarrow $ $t \in (0,t_\lambda^{e,+}(u))$ or $t\in (t_\lambda^{e,-}(u),\infty)$, $\forall \lambda \in(0,\lambda^e)$.
\end{enumerate}
	\end{lem}
	\begin{proof} The equality
				$\Lambda^{e}(tu)=\Lambda^n(tu)$  is equivalent to 
				$$
				 t^{p-q}\| u\|^p_1-\frac{(\gamma-\alpha)}{(p-\alpha)}t^{\gamma-q}\|u\|^\gamma_{L^\gamma}=q\left(t^{p-q}\frac{1}{p} \| u\|^p_1-t^{\gamma-q}\frac{(\gamma-\alpha)}{\gamma(p-\alpha)}\|u\|^\gamma_{L^\gamma}\right).
				$$
Hence,
				$$
				0=\frac{(p-q)(p-\alpha)}{p} t^{1-q}\| u\|^p_1-\frac{(\gamma-q)(\gamma-\alpha)}{\gamma}t^{\gamma-q-1}\|u\|^\gamma_{L^\gamma}=(\Lambda^{e}(tu))',
				$$
which implies \rm{$(i)$}.
		
Observe, $\mathcal{R}_\lambda ^e(tu)=\mathcal{R}_\lambda^n(tu)$ for $t>0$ if and only if 
				$$t^{p-\alpha}\| u\|^p_1+\lambda t^{q-\alpha}\|u\|^q_{L^q}-t^{\gamma-\alpha}\|u\|^\gamma_{L^\gamma}=\frac{\alpha t^{p-\alpha}}{p}\| u\|^p_1+\frac{\lambda\alpha t^{q-\alpha}}{q}-\frac{\alpha t^{\gamma-\alpha}}{\gamma}\|u\|^\gamma_{L^\gamma},$$
which is equivalent to
						\begin{align*}
					0&=\frac{(p-\alpha)}{p}t^{p-\alpha}\| u\|^p_1-\frac{\gamma-\alpha}{\gamma}t^{\gamma-\alpha}\|u\|^\gamma_{L^\gamma}-\frac{\lambda(\alpha-q)}{q}t^{q-\alpha}\|u\|^q_{L^q}\\
					&=\frac{(\alpha-q)\|u\|^q_{L^q}t^{q-\alpha}}{q}\left(\Lambda_{e}(tu)-\lambda\right).
				\end{align*}
Since $(\mathcal{R}^e_\lambda)'(tu)=0$ $\Leftrightarrow$ $\Lambda^e(tu)=0$, we get \rm{$(ii)$}.	
The proof of  \rm{$(iii)$} and \rm{$(iv)$} follow from items \rm{$(i)$},\rm{$(ii)$ and simple accounts}.	
\end{proof} 

%\begin{figure}[ht]
	%\centering
%\includegraphics[width=0.4\linewidth]{Phi.eps}\\
	%\caption{{\it Fibering function $\Phi_{\lambda,\mu}(su)$, $s\geq 0$,  $u \in W^{1,p}_0 $.}}
	%\label{fig1}
%\end{figure}
We need also 
\begin{cor}\label{muSemiCont} 
	The functionals $\lambda^e(u)$, $\lambda^n(u)$ and $\mu^{n,\pm}_\lambda(u)$ for $\lambda<\lambda^n$,  $\mu^{e,\pm}_\lambda(u)$ for $\lambda<\lambda^e$  are weakly lower semi-continuous on $ W^{1,p}_0$.  Furthermore, $t_\lambda^{n,-}(u), t_\lambda^{e,-}(u)$ are lower semi-continuous, while $t_\lambda^{n,+}(u), t_\lambda^{e,+}(u)$  upper semi-continuous on $ W^{1,p}_0$
\end{cor}
\begin{proof}
The weakly lower semi-continuity of $\lambda^e(u)$, $\lambda^n(u)$ on $ W^{1,p}_0$	follow straightforward from the definition of its (see \eqref{const0}). 

 Let us prove  as an example that  $\mu^{n,+}_\lambda(u)$ is weakly lower semi-continuous on $ W^{1,p}_0$. Assume that $\lambda<\lambda^n$.
Let $(u_m)$ be a sequence in $ W^{1,p}_0$ such that $u_m \rightharpoonup \bar{u}$ as $m \to +\infty$ weakly in $ W^{1,p}_0$ for some $\bar{u} \neq 0$.  Then $\lambda<\lambda^n<\lambda^n(\bar{u})$ and  by the weakly lower semi-continuity of $\lambda^n(u)$  we have $\lambda<\lambda^n<\lambda^n(\bar{u})\leq \liminf_{m\to +\infty}\lambda^n(u_m)$. Hence by Proposition \ref{propRN2}, there exist $t_\lambda^{n,\pm}(\bar{u})$, $t_\lambda^{n,\pm}(u_m) \in (0, +\infty)$, $m=1,2,\ldots$. Moreover, since $\Lambda^n(t\bar{u})\leq \Lambda^n(tu_m)$, for $t>0$ and for sufficiently large $m=1,2,\ldots$, we have
\begin{equation}\label{lusemicont}
	t_\lambda^{n,+}(u_m)\leq t_\lambda^{n,+}(\bar{u})<t_\lambda^{n,-}(\bar{u})\leq t_\lambda^{n,-}(u_m).  
\end{equation}
Beside this we have $\mathcal{R}_\lambda^n(t\bar{u})\leq \liminf_{m\to +\infty}\mathcal{R}_\lambda^n(tu_m)$, for all $t>0$. Hence and using the fact that $(\mathcal{R}_\lambda^n)'(t\bar{u})<0$ for $t \in (0, t_\lambda^{n,+}(\bar{u})) $ and \eqref{lusemicont} we infer
\begin{align*}
	\mu^{n,+}_\lambda(\bar{u})=\mathcal{R}_\lambda^n(t_\lambda^{n,+}(\bar{u})\bar{u})\leq & \liminf_{m\to +\infty}\mathcal{R}_\lambda^n(t_\lambda^{n,+}(u_m)\bar{u})\\
	&\leq \liminf_{m\to +\infty} \mathcal{R}_\lambda^n(t_\lambda^{n,+}(u_m)u_m)=\liminf_{m\to +\infty}\mu^{n,+}_\lambda(u_m).
\end{align*}
The proof of the last part of the corollary follows from \eqref{lusemicont}.
 \end{proof}

\section{Proof of Lemma \ref{lem1}}

First we prove
\begin{lem}\label{lem22}  Assume $1<q<\alpha<p<\gamma<  p^*$. Then,
\begin{enumerate}
	\item[${\rm (I)}$] for any $\lambda \in (0,\lambda^e)$, 
	\begin{enumerate}
	\item[{\rm $(i)$}]   \eqref{muplEn} has a minimizer $ u^{e,+}_\lambda \in W^{1,p}_0\setminus 0$ such that 
 $0<\mu^{e,+}_\lambda=\mu^{e,+}_\lambda(u^{e,+}_\lambda)$, $(\mathcal{R}^e_\lambda)''(u^{e,+}_\lambda)>0$;
	\item[{\rm $(ii)$}]    \eqref{muminusEn} has a minimizer $ u^{e,-}_\lambda \in W^{1,p}_0\setminus 0$ such that 
 $0<\mu^{e,-}_\lambda=\mu^{e,-}_\lambda(u^{e,-}_\lambda)$, $(\mathcal{R}^e_\lambda)''(u^{e,-}_\lambda)<0$;
\end{enumerate}
\item[${\rm (II) }$] for any $\lambda \in (0,\lambda^n)$,
\begin{enumerate}
\item[{\rm$(i)$}]    \eqref{muNPlus} has a minimizer $ u^{n,+}_\lambda \in W^{1,p}_0\setminus 0$  such that 
 $0<\mu^{n,+}_\lambda=\mu^{n,+}_\lambda(u^{n,+}_\lambda)$, $(\mathcal{R}^n_\lambda)''(u^{n,+}_\lambda)>0$;
\item[{\rm $(ii)$}]    \eqref{muminus} has a minimizer $ u^{n,-}_\lambda \in W^{1,p}_0\setminus 0$  such that 
 $0<\mu^{n,-}_\lambda=\mu^{n,-}_\lambda(u^{n,-}_\lambda)$, $(\mathcal{R}^n_\lambda)''(u^{n,-}_\lambda)<0$.
\end{enumerate}
\end{enumerate}
\end{lem}
		\begin{proof} The proofs of these assertions are similar. Let us prove as an example assertion \textit{(i)}, (I). 
		
	Let  $\lambda>0$. Define
the set 
$
\mathcal{Z}_\lambda:=\{u \in W^{1,p}_0\setminus 0: (\mathcal{R}_\lambda^e)'(u)=0\}. 
$  
Notice that 
\begin{align}\label{ZaitIneq}
			\mathcal{R}^e(u)=\alpha\frac{\left[\frac{(\gamma-p)}{p}\|u\|_1^p+\lambda\frac{(\gamma-q)}{q}\|u\|_{L^q}^q\right]}{(\gamma-\alpha)\|u\|_{L^{\alpha}}^\alpha}, ~\forall u \in \mathcal{Z}_\lambda.
		\end{align}
Hence	by the Sobolev inequality we derive that $\mathcal{R}^e(u)\geq \alpha\frac{(\gamma-p)}{p(\gamma-\alpha)}\|u\|_1^{p-\alpha}\to \infty$ if $u \in \mathcal{Z}_\lambda$ and $\|u\|_1 \to +\infty$. Thus, $\mathcal{R}_\lambda^e$ is coercive on $\mathcal{Z}_\lambda$, $\forall \lambda >0$. 

  Let  $(u_m)$ be a minimizing sequence of \eqref{muplEn}, i.e., 
	$
	\mu^{e,+}_\lambda(u_m)=\mathcal{R}^e_\lambda(u_m)\rightarrow \mu^{e,+}_\lambda, 
	$
where by the homogeneity of $\mu^{e,+}_\lambda(u)$ we may assume that $t_m=t_\lambda^{e,+}(u_m)=1$ for $m=1,2,\ldots$. 
The coerciveness of $\mathcal{R}_\lambda^e$  on $\mathcal{Z}_\lambda$ implies that the minimizing sequence $(u_m)$ of  \eqref{muplEn} has a weak in $W^{1,p}_0$ and strong in $L^r$, $1<r<p^*$ limit point $u^e_\lambda \in  W^{1,p}_0$.

Let us show that $u^e_\lambda \neq 0$.  
Observe, 
\begin{align}\label{lim0aB}
\mathcal{R}^e_\lambda( u_m)=\mathcal{R}^0_\lambda( u_m)-\frac{\alpha\|u_m\|^\gamma_{L^\gamma}}{ \gamma \|u_m\|^\alpha_{L^\alpha}}, ~~m=1,\ldots,
\end{align}
where 
$$
\mathcal{R}^0_\lambda(u):=\frac{\frac{1}{p}\|u\|_{1}^p+\frac{\lambda}{q}\|u\|_{L^q}^q}{\frac{1}{\alpha}\|u\|^\alpha_{L^\alpha}},~~u\in W^{1,p}_0\setminus 0.
$$
It can be shown (see e.g. \cite{DiazHerYA}) that 
\begin{equation}\label{NL}
	\min_{u \in  W^{1,p}_0\setminus 0}\min_{t>0}\mathcal{R}^0_\lambda(tu)=\mu_0>0.
\end{equation}
Denote 
$\displaystyle{a_m:= \frac{\alpha\|u_m\|^\gamma_{L^\gamma}}{ \gamma \|u_m\|^\alpha_{L^\alpha}}}$,  $m=1,2,\ldots$. 
Then 
	\begin{align*}
	\mathcal{R}^0_\lambda(u_m)= a_m \gamma \frac{\frac{1}{p}\|u_m\|_{1}^p+\frac{\lambda}{q}\|u_m\|_{L^q}^q}{\|u_m\|_{L^\gamma}^\gamma}.
\end{align*}		
We may assume  that  $a_0	:=\lim_{m\to +\infty}a_m \geq 0$. Since $\mu_0>0$, $a_0\neq 0$. Suppose, conversely to our claim, that $u^e_\lambda=0$. Then   
$$
\mathcal{R}^e_\lambda( u_m) \geq a_m \left( \gamma\frac{C}{\|u_m\|^{\gamma-p}_{L^\gamma}} - 1\right)\to +\infty,
$$
which is a contradiction. 
Thus,  $u^e_\lambda \neq 0$ and  $\mu^{e,+}_\lambda >0$.
Now by the weakly lower semi-continuity of $\mu^{e,+}_\lambda( u)$ it follows that
$\mu_\lambda^{e,+}(u^{e,+}_\lambda)\leq \mu^{e,+}_\lambda$, which obviously implies $\mu_\lambda^{e,+}=\mu_\lambda^{e,+}(u^{e,+}_\lambda)$.

\end{proof}

\begin{cor}\label{cor: ineqStr}  
%There holds
	\begin{enumerate}
		\item[{\rm $(i)$}] if $\lambda \in (0,\lambda^e)$, then  $0< \mu^{e,+}_\lambda< \mu^{e,-}_\lambda<+\infty$,
		\item[{\rm $(ii)$}] if $\lambda \in (0,\lambda^n)$, then  $0< \mu^{n,+}_\lambda< \mu^{n,-}_\lambda<+\infty$.
	\end{enumerate}
\end{cor}
\begin{proof} By  Lemma \ref{lem22}, we have 
	$$
	0<  \mu^{e,+}_\lambda\leq \mu_\lambda^{e,+}(u^{e,-}_\lambda)< \mu_\lambda^{e,-}(u^{e,-}_\lambda)=\mu^{e,-}_\lambda<+\infty.
	$$
Thus, we get  \textit{(i)}. The proof of	\textit{(ii)} is similar. 
\end{proof}
		\begin{cor}\label{mu-compare2} 
		If $\lambda\in(0,\lambda^e)$, then {\rm $(i)$} $\mu^{e,-}_\lambda<\mu^{n,-}_\lambda$; {\rm $(ii)$} $\mu^{n,+}_\lambda<\mu^{e,+}_\lambda$
						 \end{cor}
\begin{proof}
By Lemma \ref{lem22}, there exists $u^{n,-}_\lambda$ such that $\mu^{n,-}_\lambda=\mu^{n,-}_\lambda(u^{n,-}_\lambda)$. Lemma \ref{lemGEOM} entails that $\mathcal{R}^e_\lambda(t_\lambda^{e,-}(u^{n,-}_\lambda)u^{n,-}_\lambda)$ = $\mathcal{R}^n_\lambda(t_\lambda^{e,-}(u^{n,-}_\lambda)u^{n,-}_\lambda)$ and the function $t\mapsto \mathcal{R}^n_\lambda(tu^{n,-}_\lambda)$ is decreasing on the interval  $(t_\lambda^{n,-}(u^{n,-}_\lambda),t_\lambda^{e,-}(u^{n,-}_\lambda))$.
Hence, 
				$$\mu_\lambda^{e,-}\leq\mathcal{R}^e_\lambda(t_\lambda^{e,-}(u^{n,-}_\lambda)u^{n,-}_\lambda)=\mathcal{R}^n_\lambda(t_\lambda^{e,-}(u^{n,-}_\lambda)u^{n,-}_\lambda)<\mathcal{R}^n_\lambda(t_\lambda^{n,-}(u^{n,-}_\lambda)u^{n,-}_\lambda)=\mu^{n,-}_\lambda,
				$$
and we get 		\textit{(i)}.		 The proof of \textit{(ii)} is similar. 
	\end{proof}
			
	\medskip
			
\noindent {\textit{The proof of Lemma \ref{lem1}}:}
		 follows from Corollaries \ref{cor: ineqStr}, \ref{mu-compare2}.

 \medskip

\begin{cor}\label{muCor} 
		
		\begin{enumerate}
		\item[{\rm $(i)$}] The minimizer $ u^{e,+}_\lambda $ of \eqref{muplEn} (perhaps, after a scaling) is a non-negative critical point of $\Phi_{\lambda,\mu^{e,+}_\lambda}$, moreover   $\Phi_{\lambda,\mu^{e,+}_\lambda}(u^{e,+}_\lambda)=0$ and  $\Phi_{\lambda,\mu^{e,+}_\lambda}''(u^{e,+}_\lambda)>0$.
			\item[{\rm $(ii)$}] The minimizer $ u^{e,-}_\lambda$ of \eqref{muminusEn} (perhaps, after a scaling) is a non-negative critical point of $\Phi_{\lambda,\mu^{e,-}_\lambda}$ moreover  $\Phi_{\lambda,\mu^{e,-}_\lambda}(u^{e,-}_\lambda)=0$ and $\Phi_{\lambda,\mu^{e,-}_\lambda}''(u^{e,-}_\lambda)<0$.

		\end{enumerate}
\end{cor}
\begin{proof} \textit{(i)} Let $u^{e,+}_\lambda$ be a minimizer of \eqref{muplEn}. Then $u^{e,+}_\lambda$ is also a minimizer of \eqref{mumue} with $t_\lambda^{e,+}(u^{e,+}_\lambda)=1$. Hence,
$0=D\mu^{e,+}_\lambda(u^{e,+}_\lambda)=D\mathcal{R}^e_\lambda(u^{e,+}_\lambda)$,  
and consequently,  $D\Phi_{\lambda,\mu^{e,+}_\lambda}(u^{e,+}_\lambda)=0$.
 Moreover,  $\mathcal{R}_\lambda^e(u^{e,+}_\lambda)=\mu^{e,+}_\lambda$, $(\mathcal{R}_\lambda^e)''(u^{e,+}_\lambda)>0$ yield $\Phi_{\lambda,\mu^{e,+}_\lambda}(u^{e,+}_\lambda)=0$ and  $\Phi_{\lambda,\mu^{e,+}_\lambda}''(u^{e,+}_\lambda)>0$, respectively. Since $\mu^{e,+}_\lambda(|u^{e,+}_\lambda|) =\mu^{e,+}_\lambda(u^{e,+}_\lambda)$  one may assume that $u^{e,+}_\lambda\geq 0$. The proof of $(ii)$ is similar.
	
\end{proof}

\section{Proof of $(1^o)$, Theorem \ref{thm1}}

%Let us show that $\Phi_{\lambda,\mu}(u^1_{\lambda,\mu}) \to 0$ as $\mu \to +\infty$. 

We obtain the solution $u^2_{\lambda,\mu}$ using the following Nehari minimization problem
\begin{equation} \label{GminUnst0}
	\hat{\Phi}^2_{\lambda,\mu}=\min \{\Phi_{\lambda,\mu}(u):  ~ u \in \mathcal{R}\mathcal{N}_{\lambda,\mu}^2\},
\end{equation}
where 
$
\mathcal{R}\mathcal{N}_{\lambda,\mu}^2:=\left\{u\in \mathcal{N}_{\lambda,\mu}:~ (\mathcal{R}_\lambda^{n})'(u)\geq 0 \right\}.
$ 
Observe,  $\mathcal{R}\mathcal{N}_{\lambda,\mu}^2\neq \emptyset$, for all $\lambda \in (0,\lambda^e)$ and  $\mu \in (\mu^{n,+}_\lambda, \mu^{n,-}_\lambda)$. Indeed, if $\mu \in (\mu^{n,+}_\lambda, \mu^{n,-}_\lambda)$, then  $\mu<\mu^{n,-}_\lambda<\mu^{n,-}_\lambda(u)$,  $\forall u \in W^{1,p}_0\setminus 0$, and therefore, there exists $\tilde{u} \in W^{1,p}_0\setminus 0$ such that $\mu^{n,+}_\lambda<\mu^{n,+}_\lambda(\tilde{u})<\mu<\mu^{n,-}_\lambda(\tilde{u})$.  Lemma \ref{lemGEOM} implies that there exists $s^2_{\lambda, \mu}(\tilde{u}) \in  (t^{n,+}_\lambda(\tilde{u}), t^{n,-}_\lambda(\tilde{u}))$, and thus, $s^2_{\lambda, \mu}(\tilde{u})\tilde{u} \in \mathcal{R} \mathcal{N}_{\lambda,\mu}^2$. 

Furthermore, the assumption $\lambda \in (0,\lambda^e)$, $\mu \in (\mu^{n,+}_\lambda, \mu^{n,-}_\lambda)$ implies that
\begin{equation}\label{main:ineq}
	\mu^{n,+}_\lambda(u)\leq \mu<\mu^{n,-}_\lambda(u), ~~\forall u \in \mathcal{R}\mathcal{N}_{\lambda,\mu}^2.
\end{equation}
Indeed, $\mu<\mu^{n,-}_\lambda\leq \mu^{n,-}_\lambda(u)$  for any $u \in  W^{1,p}_0\setminus 0$, whereas the conditions $\mathcal{R}_\lambda^{n}(u)=\mu$, $(\mathcal{R}_\lambda^{n})'(u)\geq 0$  for $u \in \mathcal{R}\mathcal{N}_{\lambda,\mu}^2$ yield $\mu^{n,+}_\lambda(u)\leq \mu$.

\begin{lem}\label{lem:22_1}
	Let $\lambda \in (0,\lambda^e)$ and  $\mu \in (\mu^{n,+}_\lambda, \mu^{n,-}_\lambda)$. There exists a minimizer $\bar{u}^2_\mu \in W^{1,p}_0\setminus 0$ of \eqref{GminUnst0} and
		\begin{equation}\label{eq:ineq+}
		\begin{cases}
		& \hat{\Phi}^2_{\lambda,\mu}=\Phi_{\lambda,\mu}(\bar{u}^2_\mu) >0~~\mbox{if}~~  ~\mu \in (\mu^{n,+}_\lambda, \mu^{e,+}_\lambda),\\
	& \hat{\Phi}^2_{\lambda,\mu}=\Phi_{\lambda,\mu}(\bar{u}^2_\mu) =0~~\mbox{if}~~  ~\mu= \mu^{e,+}_\lambda,\\
			&\hat{\Phi}^2_{\lambda,\mu}=\Phi_{\lambda,\mu}(\bar{u}^2_\mu)<0~~\mbox{if}~~   ~\mu \in (\mu^{e,+}_\lambda, \mu^{n,-}_\lambda).	
				\end{cases}
\end{equation}
		
\end{lem}
\begin{proof}
Note that by Lemma \ref{lem22}, the minimum of \eqref{GminUnst0} for $\mu =\mu^{e,+}_\lambda$ attains (perhaps, after a scaling) at the solution  $ u^{e,+}_\lambda \in W^{1,p}_0\setminus 0$ of \eqref{muplEn} and  $\hat{\Phi}^2_{\lambda,\mu}|_{\mu =\mu^{e,+}_\lambda}=0$. This easily implies the proof of the lemma in the case  $\mu=\mu^{e,+}_\lambda$.

	Assume that $\mu \in (\mu^{n,+}_\lambda,\mu^{n,-}_\lambda)$.	Let $u_m$ be a minimizing sequence of \eqref{GminUnst0}. The coerciveness of $\Phi_{\lambda,\mu}$  implies that the sequence $(u_m)$ is bounded in $W^{1,p}_0$ and thus,  up to a
	subsequence, 
	$$
	u_m \to \bar{u}^2_\mu~~\mbox{strongly in}~~L^r,~~\mbox{and weakly in} ~W^{1,p}_0,
	$$
	for some $\bar{u}^2_\mu \in W^{1,p}_0$, where $r \in(1,p^*)$. Moreover,  $\Phi_{\lambda,\mu}(\bar{u}^2_\mu)\leq \liminf_{m\to \infty}\Phi_{\lambda, \mu}(u_m)=\hat{\Phi}^2_{\lambda,\mu}$.

If $\mu \in (\mu^{e,+}_\lambda, \mu^{n,-}_\lambda)$, then there exists $\tilde{u}\in \mathcal{R}\mathcal{N}_{\lambda,\mu}^2$ such that $1=s^2_{\lambda, \mu}(\tilde{u}) \in  (t^{e,+}_\lambda(\tilde{u}), t^{e,-}_\lambda(\tilde{u}))$. Hence by \textit{(iv)}, Lemma \ref{lemGEOM},   $\mu=\mathcal{R}_\lambda^n(s^2_{\lambda,\mu}(\tilde{u})\tilde{u})>\mathcal{R}_\lambda^e(s^2_{\lambda,\mu}(\tilde{u})\tilde{u})\equiv \mathcal{R}_\lambda^e(\tilde{u})$, which implies   $0> \Phi_{\lambda, \mu}(\tilde{u})\geq\hat{\Phi}^2_{\lambda,\mu}$.  Hence, $\Phi_{\lambda,\mu}(\bar{u}^2_\mu)<0$, and consequently,   $\bar{u}^2_\mu\neq 0$ for $\mu \in (\mu^{e,+}_\lambda, \mu^{n,-}_\lambda)$.

If  $\mu \in (\mu^{n,+}_\lambda, \mu^{e,+}_\lambda)$, then $\mu<\mu^{e,+}_\lambda<\mu^{e,+}_\lambda(u)=\mathcal{R}_\lambda^e(t^{e,+}_\lambda(u)u) $ for any $u \in \mathcal{R}\mathcal{N}_{\lambda,\mu}^2$, and therefore, $1=s^2_{\lambda, \mu}(u) \in  (0,t^{e,+}_\lambda(u))$.  Hence by \textit{(iv)}, Lemma \ref{lemGEOM}, 
 $\mathcal{R}_\lambda^e(u)\equiv\mathcal{R}_\lambda^e(s^2_{\lambda,\mu}(u)u)>\mu$,  and consequently, $ \Phi_{\lambda, \mu}(u) >0$,   $\forall \mu \in (\mu^{n,+}_\lambda, \mu^{e,+}_\lambda)$, $\forall u \in \mathcal{R}\mathcal{N}_{\lambda,\mu}^2$. 
Let us show that $\bar{u}^2_\mu\neq 0$ for $\mu \in (\mu^{n,+}_\lambda, \mu^{e,+}_\lambda)$. To obtain a contradiction, suppose that $\bar{u}^2_\mu= 0$. Observe that if  $\mu \in (\mu^{n,+}_\lambda, \mu^{e,+}_\lambda)$, then $t^{n,+}_\lambda(u_m)<s^2_{\lambda,\mu}(u_m)=1$, $m=1,2,\ldots$. Consequently, $t^{n,+}_\lambda(u_m)u_m \to 0$ as $m\to +\infty$ strongly in $L^r$,  $r \in(1,p^*)$ and weakly in $W^{1,p}_0$.
Analysis similar to that in the proof of Lemma \ref{lem22}  shows that this implies $\mathcal{R}_\lambda^n(t^{n,+}_\lambda(u_m)u_m) \to +\infty$. However, this contradicts $\mathcal{R}_\lambda^n(t^{n,+}_\lambda(u_m)u_m)\leq \mathcal{R}_\lambda^n(u_m)=\mu$, $m=1,2,\ldots$, and thus, $\bar{u}^2_\mu\neq 0$ for $\mu \in (\mu^{n,+}_\lambda, \mu^{e,+}_\lambda)$.

Let $\mu \in (\mu^{n,+}_\lambda, \mu^{e,+}_\lambda)$ and suppose, contrary to our claim, that $\Phi_{\lambda,\mu}(u_m)\to 0$ as $m\to +\infty$. Then 
\begin{equation}
\Phi_{\lambda,\mu}(u_m)=\| u_m\|^q_{L^q}\left(\mathcal{R}^e_\lambda(u_m)-\mu\right)\to 0.
\end{equation}
Since $\Phi_{\lambda,\mu}(u_m)>0$, $m=1,2,\ldots$ and $\bar{u}^2_\mu\neq 0$, this implies $\mathcal{R}^e_\lambda(u_m) \downarrow \mu$ as $m\to +\infty$. Since $\mu<\mu^{e,+}_\lambda\leq \mu^{e,+}_\lambda(u_m)=\mathcal{R}_\lambda^e(t^{e,+}_\lambda(u_m)u_m) $, $1=s^2_{\lambda, \mu}(u_m) \in  (0,t^{e,+}_\lambda(u_m))$, and therefore, $\mathcal{R}^e_\lambda(u_m)>\mathcal{R}_\lambda^e(t^{e,+}_\lambda(u_m)u_m)=\mu^{e,+}_\lambda(u_m)>\mu$, $m=1,2,\ldots$. Consequently,   $\lim_{m\to +\infty}\mu^{e,+}_\lambda(u_m)= \mu <\mu^{e,+}_\lambda$, which is a contradiction since $\mu^{e,+}_\lambda(u_m)\geq \mu^{e,+}_\lambda$, $m=1,2,\ldots$. Thus, $\hat{\Phi}^2_{\lambda,\mu} >0~~\mbox{if}~~  ~\mu \in (\mu^{n,+}_\lambda, \mu^{e,+}_\lambda)$, and we have proved \eqref{eq:ineq+}.

Let us show that $\bar{u}^2_\mu$ is a minimizer of \eqref{GminUnst0}, i.e., $\bar{u}^2_\mu \in \mathcal{R}\mathcal{N}_{\lambda,\mu}^2$ and $\Phi_{\lambda, \mu}(\bar{u}^2_\mu )= \hat{\Phi}^2_{\lambda, \mu}$. 
By \eqref{main:ineq}, $\mu^{n,+}_\lambda(\bar{u}^2_\mu)\leq \mu<\mu^{n,-}_\lambda(\bar{u}^2_\mu)$, and therefore
 $\exists s^2_{\lambda, \mu}(\bar{u}^2_\mu)\in (s^1_{\lambda, \mu}(\bar{u}^2_\mu), s^3_{\lambda, \mu}(\bar{u}^2_\mu))$ such that  
	\begin{equation*}\label{ineqPHI}
	\begin{aligned}
		&\mu=\mathcal{R}^n_\lambda(s^2_{\lambda, \mu}(\bar{u}^2_\mu)\bar{u}^2_\mu)\leq \liminf_{m\to \infty}\mathcal{R}^n_\lambda(s^2_{\lambda, \mu}(\bar{u}^2_\mu)u_m),\\
			&0<(\mathcal{R}^n_\lambda)'(s^2_{\lambda, \mu}(\bar{u}^2_\mu)\bar{u}^2_\mu)\leq  \liminf_{m\to \infty}(\mathcal{R}^n_\lambda)'(s^2_{\lambda, \mu}(\bar{u}^2_\mu)u_m).
				\end{aligned}
	\end{equation*}
	This means that $1=s^2_{\lambda,\mu}(u_m)\leq s^2_{\lambda, \mu}(\bar{u}^2_\mu) <s^3_{\lambda, \mu}(u_m)$, $m=1,\ldots $. Hence by
	$$
	\mathcal{R}^n_\lambda(\bar{u}^2_\mu)\leq \liminf_{m\to \infty}\mathcal{R}^n_\lambda(u_m) = \mu,
	$$
	we obtain $s^1_{\lambda, \mu}(\bar{u}^2_\mu) \leq 1\leq s^2_{\lambda, \mu}(\bar{u}^2_\mu) $. Since $\Phi'_{\lambda, \mu}(s \bar{u}^2_\mu)<0$,  for any $s \in (s^1_{\lambda, \mu}(\bar{u}^2_\mu),s^2_{\lambda, \mu}(\bar{u}^2_\mu))$, we derive
	$$
	\Phi_{\lambda,\mu}(s^2_{\lambda, \mu}(\bar{u}^2_\mu)\bar{u}^2_\mu)\leq \Phi_{\lambda,\mu}(\bar{u}^2_\mu)\leq \liminf_{m\to \infty}\Phi_{\lambda, \mu}(u_m)=\hat{\Phi}^2_{\lambda,\mu},
	$$
	which yields that  $s^2_{\lambda, \mu}(\bar{u}^2_\mu)\bar{u}^2_\mu=\bar{u}^2_\mu$ is a minimizer of  \eqref{GminUnst0} and thus, $u_m \to \bar{u}^2_\mu$ strongly in $W^{1,p}_0$.

\end{proof}
Consider the following subset of $\mathcal{M}_{\lambda, \mu}$
$$
\mathcal{M}^2_{\lambda, \mu}:=\{u \in \mathcal{R}\mathcal{N}_{\lambda,\mu}^2 :~~ \hat{\Phi}^2_{\lambda,\mu}=\Phi_{\lambda,\mu}(u) \}.
$$
Lemma \ref{lem:22_1} yields that  $\mathcal{M}^2_{\lambda, \mu} \neq \emptyset$ for any $\lambda \in (0,\lambda^e)$ and  $\mu \in (\mu^{n,+}_\lambda, \mu^{n,-}_\lambda)$. Note that the minimizer $\bar{u}^2_\mu\in \mathcal{M}^2_{\lambda, \mu}$ of \eqref{GminUnst0}  does not necessarily provide a solution of \eqref{p}. By Lemma \ref{lemNL} and \eqref{RandPhi},  $\bar{u}^2_\mu\in \mathcal{M}^2_{\lambda, \mu}$  corresponds to a solution of \eqref{p}, if the strict inequality $(\mathcal{R}_\lambda^{n})'(\bar{u}^2_\mu)>0$ holds. Note that by	Proposition \ref{propRN2}, 
$$
(\mathcal{R}_\lambda^{n})'(\bar{u}^2_\mu)>0~ \Leftrightarrow~\mu^{n,+}_\lambda(\bar{u}^2_\mu)<\mu<\mu^{n,-}_\lambda(\bar{u}^2_\mu).
$$
By \eqref{main:ineq}, if $\lambda \in (0,\lambda^e)$ and  $\mu \in (\mu^{n,+}_\lambda, \mu^{n,-}_\lambda)$, then  $\mu<\mu^{n,-}_\lambda(\bar{u}^2_\mu)$ for any $\bar{u}^2_\mu\in \mathcal{M}^2_{\lambda, \mu}$. Thus, to obtain that  $\bar{u}^2_\mu\in \mathcal{M}^2_{\lambda, \mu}$ is a weak solution of \eqref{p} it is sufficient to show that $\mu^{n,+}_\lambda(\bar{u}^2_\mu)<\mu$.

\begin{cor}\label{cor:mu1}
	Let $\lambda \in (0,\lambda^e)$. If $\mu \in [\mu^{e,+}_\lambda, \mu^{n,-}_\lambda)$, then  $\bar{u}^2_\mu\in \mathcal{M}^2_{\lambda, \mu}$ is a weak solution of \eqref{p}.
\end{cor}
\begin{proof}
 Indeed, the inequality $\Phi_{\lambda, \mu}(\bar{u}^2_\mu )=\hat{\Phi}^2_{\lambda, \mu}\leq 0$ implies  $\mu^{n,+}_\lambda(\bar{u}^2_\mu)<\mu^{e,+}_\lambda(\bar{u}^2_\mu)\leq \mathcal{R}^e_\lambda(\bar{u}^2_\mu)\leq \mu$, and thus, $\mu^{n,+}_\lambda(\bar{u}^2_\mu)<\mu$. 
\end{proof}
\begin{lem}\label{lemCC3}
Let $\lambda \in (0,\lambda^e)$.  There exists $\tilde{\mu}_\lambda^{*}\in (\mu^{n,+}_\lambda, \mu^{e,+}_\lambda)$ such that 	if $\mu \in (\tilde{\mu}_\lambda^{*},\mu^{e,+}_\lambda)$, then  $\bar{u}^2_{\mu} \in \mathcal{M}^2_{\lambda, \mu}$ weakly satisfies \eqref{p}, moreover  $\Phi_{\lambda, \mu}(\bar{u}^2_\mu )= \hat{\Phi}^2_{\lambda, \mu}>0$,  $\Phi_{\lambda, \mu}''(\bar{u}^2_\mu)>0$.
\end{lem}
\begin{proof} By the above, it is sufficient to show that there exists $\tilde{\mu}_\lambda^{*} \in (\mu^{n,+}_\lambda,  \mu^{e,+}_\lambda)$ such that  $\mu^{n,+}_\lambda(u)<\mu$, $\forall u \in \mathcal{M}^2_{\lambda, \mu}$, $\forall \mu \in (\tilde{\mu}_\lambda^{*},  \mu^{e,+}_\lambda)$. 
Suppose this is false. Then there exist sequences $\mu_m \in (\mu^{n,+}_\lambda,  \mu^{e,+}_\lambda)$ and $u_m \in \mathcal{M}^2_{\lambda, \mu}$, $m=1,2,\ldots$ such that $\mu_m \to  \mu^{e,+}_\lambda$ as $m\to +\infty$ and $\mu_m= \mu^{n,+}_\lambda(u_m)$, $\forall m=1,2,\ldots$. By Proposition \ref{ContPh},  up to a subsequence, $u_m \to u_0$  strongly in $W^{1,p}_0$ as $m\to +\infty$ for some $u_0 \in \mathcal{M}^2_{\lambda,  \mu^{e,+}_\lambda}$. Hence $\mu^{n,+}_\lambda(u_0)= \mu^{e,+}_\lambda$, and consequently $u_0 \in \mathcal{R}\mathcal{N}_{\lambda, \mu^{e,+}_\lambda}^2$. Furthermore, Corollary \ref{cor:contPHI} implies that $\Phi_{\lambda, \mu^{e,+}_\lambda}(u_0)= \lim_{m\to +\infty}\Phi_{\lambda,\mu_m}(u_m)=\hat{\Phi}_{\lambda,  \mu^{e,+}_\lambda}^2$. Hence 
$u_0 \in \mathcal{M}^2_{\lambda,  \mu^{e,+}_\lambda}$ and $\mu^{n,+}_\lambda(u_0)= \mu^{e,+}_\lambda(u_0)=\mu^{e,+}_\lambda$, which contradicts ($2^o$), Lemma \ref{lem1}. 
\end{proof}

\textit{Conclusion of the proof of $(1^o)$, Theorem \ref{thm1}. 
}

Let $0<\lambda<\lambda^e$. Introduce,
\begin{equation}
\hat{\mu}_\lambda^{*}:=\sup\{\mu \in (\mu^{n,+}_\lambda, \mu^{n,-}_\lambda): ~	\mu= \mu^{n,+}_\lambda(u), ~ \exists u \in \mathcal{M}^2_{\lambda, \mu}\}.
\end{equation}
Corollary \ref{cor:mu1} and Lemma \ref{lemCC3} imply that $ \mu^{n,+}_\lambda\leq \hat{\mu}_\lambda^{*} <\mu^{e,+}_\lambda$. Hence for any $\mu \in (\hat{\mu}_\lambda^{*},\mu^{n,-}_\lambda)$, there holds  $\mu^{n,+}_\lambda(u)<\mu^{e,+}_\lambda$,	$\forall u \in \mathcal{M}^2_{\lambda, \mu^{e,+}_\lambda}$, and therefore, each $u^2_{\lambda, \mu}:=\bar{u}^2_{\mu} \in \mathcal{M}^2_{\lambda, \mu}$ is a weak solution of \eqref{p}. By the above,    $\Phi_{\lambda, \mu}''(u_{\lambda,\mu}^2 )>0$, $\Phi_{\lambda,\mu}(u^2_{\lambda, \mu})>0$, for $\mu\in (\hat{\mu}_\lambda^{*},\mu^{e,+}_\lambda)$, $\Phi_{\lambda,\mu^{e,+}_\lambda}(u^2_{\lambda, \mu^{e,+}_\lambda})=0$, and $\Phi_{\lambda,\mu}(u^2_{\lambda, \mu})<0$ for $\mu\in (\mu^{e,+}_\lambda,\mu^{n,-}_\lambda)$.
From the above, $u^2_{\lambda, \mu}$ is a local minimizer of $\Phi_{\lambda,\mu}(u)$ in the Nehari manifold
$\mathcal{N}_{\lambda,\mu}$. This by $\Phi_{\lambda, \mu}''(u_{\lambda,\mu}^2 )>0$ implies that  $u^2_{\lambda, \mu}$ is a local minimizer of $\Phi_{\lambda,\mu}(u)$ in $W^{1,p}_0$. Thus, $u^2_{\lambda,\mu}$ is a linearly stable solution. It is obvious that $u^2_{\lambda, \mu}$ is a ground state.

Note that $\Phi_{\lambda,\mu}(|u_{\lambda,\mu}^2 |) =\Phi_{\lambda,\mu}(u_{\lambda,\mu}^2 )$ and $|u_{\lambda,\mu}^2 | \in \mathcal{R}\mathcal{N}_{\lambda,\mu}^2$. Hence  one may assume that $u_{\lambda,\mu}^2 \geq 0$.  
The bootstrap argument	and the Sobolev embedding theorem
	yield that $u^2_{\lambda, \mu} \in L^\infty$. Then $C^{1,\kappa}$-regularity results of DiBenedetto \cite{DiBened} and Tolksdorf \cite{Tolksdorf} (interior regularity) combined with  Lieberman \cite{Lieber} (regularity up to
	the boundary) yield $u^2_{\lambda, \mu} \in C^{1,\kappa}(\overline{\Omega})$ for  $\kappa \in (0,1)$. Furthermore, since $p<\gamma$, the Harnack inequality due to Trudinger \cite{trudin} implies that $u^2_{\lambda, \mu}>0$ in $\Omega$.

	From Corollary \ref{cor:contPHI} it follows that the function  $(\hat{\mu}_\lambda^{*},\mu^{n,-}_\lambda) \ni \mu \mapsto \Phi_{\lambda,\mu}(u^2_{\lambda, \mu})$ is continuous and monotone decreasing.

\medskip

\section{Proof of $(2^o)$, Theorem \ref{thm1}}
\noindent

Let $\lambda \in (0,\lambda^n)$ and $\mu \in (-\infty, +\infty)$. Consider
\begin{equation} \label{GminUnst1}
	\hat{\Phi}^3_{\lambda,\mu}=\min \{\Phi_{\lambda,\mu}(u):~ u \in \mathcal{R}\mathcal{N}_{\lambda,\mu}^{3}\},
\end{equation}
where 
$$
\mathcal{R}\mathcal{N}_{\lambda,\mu}^{3}:=\{u\in \mathcal{N}_{\lambda,\mu}: ~ (\mathcal{R}^n_\lambda)'(u)\leq 0,~(\Lambda^n)'(u)<0\}.
$$
Notice that
$
\mathcal{R}\mathcal{N}_{\lambda,\mu}^{3}=\{u\in \mathcal{N}_{\lambda,\mu}: ~ (\mathcal{R}^n_\lambda)'(u)\leq 0,~t^n(u)<1\}$, where $t^n(u)$ is defined by \eqref{17}.

\begin{lem}\label{lem:09_10}
	Let $\lambda \in (0,\lambda^n)$ and  $\mu \in (-\infty, +\infty)$. Then there exists a minimizer $\bar{u}^3_{\lambda,\mu} \in \mathcal{R}\mathcal{N}_{\lambda,\mu}^{3}$ of \eqref{GminUnst1} and  $\Phi_{\lambda,\mu}''(\bar{u}^3_{\lambda,\mu}) \leq 0$.
		\end{lem}
\begin{proof}
Since $\sup_{u \in W_0^{1,p}\setminus 0}  \mu^{n,-}_\lambda(u)=+\infty$, 
 one can find $u \in W^{1,p}_0\setminus 0$ for any  $\mu \in (-\infty, +\infty)$ such that $\mu<\mu^{n,-}_\lambda(u)$,  and therefore, there exists   $s^3_{\lambda,\mu}(u)>t^n(u)$. Hence $s^3_{\lambda,\mu}(u)u \in \mathcal{R}\mathcal{N}_{\lambda,\mu}^{3}$ and therefore,  
$\mathcal{R}\mathcal{N}_{\lambda,\mu}^{3}\neq \emptyset $ for any  $\lambda \in (0,\lambda^n)$, $\mu \in (-\infty, +\infty)$. 

Let $(u_m)$ be a minimizing sequence of \eqref{GminUnst1}. Similar to the proof of $(1^o)$, Theorem \ref{thm1} one can deduce that there exists a subsequence, which we again denote   by $(u_m)$,  and a limit point $\bar{u}^3_\mu$ such that $u_m \to \bar{u}^3_\mu$ strongly in $L^r$,  $r \in(1,p^*)$ and weakly in $W^{1,p}_0$.
Observe that if $\bar{u}^3_\mu= 0$, then by \eqref{17} we obtain a contradiction 
$$
1>(t^{n}(u_m))^{(\gamma-p)}=C_n\frac{\| u_m\|_1^p}{\|u_m\|_{L^\gamma}^\gamma} \geq c \frac{1}{\| u_m\|_{L^\gamma}^{\gamma-p} }\to +\infty~~\mbox{as}~~ m\to +\infty,
$$ 
where $c \in (0,+\infty)$ does not depend on $m$. Thus $\bar{u}^3_\mu\neq  0$, and therefore, there exists $s^3_{\lambda,\mu}(\bar{u}^3_\mu)>t^n(\bar{u}^3_\mu)$ so that $s^3_{\lambda,\mu}(\bar{u}^3_\mu) \bar{u}^3_\mu \in \mathcal{R}\mathcal{N}_{\lambda,\mu}^{3}$. 

By Corollary \ref{muSemiCont}, we have $t_\lambda^{n,+}(u_m)\leq t_\lambda^{n,+}(\bar{u})<t_\lambda^{n,-}(\bar{u})\leq t_\lambda^{n,-}(u_m)$ for sufficiently large $m$. From this and since
 $\mathcal{R}^n_\lambda(t\bar{u}^3_\mu) \leq \liminf_{m\to +\infty}\mathcal{R}^n_\lambda(t u_m)$, for any $t>0$, it follows that for sufficiently large $m$ there holds  $s^3_{\lambda,\mu}(\bar{u}^3_\mu)\leq s^3_{\lambda,\mu}(u_m)=1$, and if $ s^2_{\lambda,\mu}(u_m)$ exists,  $s^2_{\lambda,\mu}(u_m)<s^2_{\lambda,\mu}(\bar{u}^3_\mu)\leq s^3_{\lambda,\mu}(\bar{u}^3_\mu)$. Hence by the weak lower semi-continuity of $\Phi_{\lambda,\mu}(u)$ we have
$$
\Phi_{\lambda,\mu}(s^3_{\lambda,\mu}(\bar{u}^3_\mu)\bar{u}^3_\mu) \leq \liminf_{m\to +\infty}\Phi_{\lambda,\mu}(s^3_{\lambda,\mu}(\bar{u}^3_\mu)u_m)\leq\liminf_{m\to +\infty}\Phi_{\lambda,\mu}(u_m)=\hat{\Phi}^3_{\lambda,\mu}, 
$$ 
which implies that $s^3_{\lambda,\mu}(\bar{u}^3_\mu) \bar{u}^3_\mu$ is a minimizer, and consequently,  $s^3_{\lambda,\mu}(\bar{u}^3_\mu)=1$ and $u_m \to \bar{u}^3_{\lambda,\mu}$ strongly in $W^{1,p}_0$.   Since $\lambda<\lambda^n$, we have $(\Lambda^n)'(\bar{u}^3_\mu)<0$, and therefore, $(\mathcal{R}^n_\lambda)'(\bar{u}^3_\mu)\leq 0$.
\end{proof}

Assume that $\mu \in (-\infty, \mu^{n,-}_\lambda)$. Then $\mu<\mu^{n,-}_\lambda<\mu^{n,-}_\lambda(\bar{u}^3_\mu)$,  and therefore, $(\mathcal{R}^n_\lambda)'(\bar{u}^3_\mu)<0$. This by Lemma \ref{lemNL} and \eqref{RandPhi} implies that $u^3_{\lambda, \mu}:=\bar{u}^3_\mu$ is a weak solution of \eqref{p}. Moreover, since $(\mathcal{R}^n_\lambda)'(u^3_{\lambda, \mu})<0$, we have $\Phi''_{\lambda,\mu 	}(u^3_{\lambda, \mu})<0$, and therefore, $u^3_{\lambda, \mu}$ is a linearly unstable solution. 
Analysis similar to that in the proof of $(1^o)$, Theorem \ref{thm1} shows that  $u^3_{\lambda,\mu} \in C^{1,\kappa}(\overline{\Omega})$ for  $\kappa \in (0,1)$ and $u^3_{\lambda,\mu}>0$. As \eqref{eq:ineq+} it can be shown that $\Phi_{\lambda,\mu}(u^3_{\lambda, \mu})>0$ if $\mu \in (-\infty, \mu^{e,-}_\lambda)$, $\Phi_{\lambda,\mu^{e,-}_\lambda}(u^3_{\lambda, \mu^{e,-}_\lambda})=0$,  and  $\Phi_{\lambda,\mu}(u^3_{\lambda, \mu})<0$ if $\mu \in (\mu^{e,-}_\lambda,\mu^{n,-}_\lambda)$. 
Corollary \ref{cor:contPHI} implies that the function  $(-\infty, \mu^{n,-}_\lambda) \ni \mu \mapsto \Phi_{\lambda,\mu}(u^3_{\lambda, \mu})$ is continuous and monotone decreasing. 

Let us show \textit{(iv)}. From the monotonicity of  $\Phi_{\lambda,\mu}(u^3_{\lambda, \mu})$ it follows
 $\Phi_{\lambda,\mu}(u^3_{\lambda, \mu}) \to C$ as $\mu \to -\infty$ for some $C\in (0,+\infty]$.  Suppose, contrary to our claim, that $C<+\infty$.
  Since $u_{\lambda,\mu}^3\in \mathcal{N}_{\lambda,\mu}$, 
				\begin{align}
			\label{limlambda}
						\frac{C}{2}<&\Phi_{\lambda,\mu}(u_{\lambda,\mu}^3)=\\
						&\frac{\gamma-p}{\gamma p}\|u_{\lambda,\mu}^3\|_1^p+\lambda\frac{\gamma-q}{\gamma q}\|u_{\lambda,\mu}^3\|_{L^q}^q-\mu\frac{\gamma-\alpha}{\gamma \alpha}\|u_{\lambda,\mu}^3\|_{L^\alpha}^\alpha<\frac{3C}{2}, \nonumber
\end{align}
				for  sufficiently large $|\mu|$. This 
		implies that $\|u_{\lambda,\mu}^3\|_{L^\alpha}^\alpha \to 0$ as $\mu \to -\infty$, and 
 $(u_{\lambda,\mu}^3)$ is bounded. Thus, there exists a subsequence  $\mu_j \to -\infty$ such that $u_{\lambda,\mu_j}^3\rightharpoonup \overline{u}$ in $W_0^{1,p}$ as $j\to +\infty$ for some  $\overline{u}\in W_0^{1,p}$. Since $\|u_{\lambda,\mu}^3\|_{L^\alpha}^\alpha \to 0$ as $\mu \to -\infty$, $\|u_{\lambda,\mu_j}^3\|_{L^q}\to 0$. Hence  passing to the limit in $\Phi_{\lambda,\mu}'(u_{\lambda,\mu_j}^3)=0$ we obtain 
 $
\lim_{j\to +\infty}\|u_{\lambda, \mu_j}^{3}\|_1^p=0
$.  
 This and \eqref{limlambda} yield  $0<{C}/{2}\leq \lim_{\mu_j \to \-\infty}\Phi_{\lambda,\mu}(u_{\lambda,\mu_j}^3)=0$, which is a contradiction. Thus $\Phi_{\lambda,\mu}(u^3_{\lambda, \mu}) \to +\infty$ and $\|u_{\lambda,\mu}^3\|_1 \to +\infty$ as $\mu \to -\infty$.  

Let us show that $u^3_{\lambda,\mu}$ is a ground state of \eqref{p} if $\mu \in (-\infty, \mu^{n,+}_\lambda)$.  By Proposition \ref{propRN2} and \eqref{mupl},   if $\mu \in (-\infty, \mu^{n,+}_\lambda)$, then for any $u \in W^{1,p}\setminus 0$, the fibering function $\Phi_{\lambda,\mu}(su)$ has only critical point $s^3_{\lambda, \mu}(u)>0$. Hence if $\mu \in (-\infty, \mu^{n,+}_\lambda)$, then $\mathcal{R}\mathcal{N}_{\lambda,\mu}^{3}= \mathcal{N}_{\lambda,\mu}$  and 
$\inf_{u \in W^{1,p}_0\setminus 0}\max_{s>0}\Phi_{\lambda,\mu}(su)=\hat{\Phi}^3_{\lambda,\mu}$ and we obtain the desired. This concludes the proof of $(2^o)$, Theorem \ref{thm1}.

\medskip

\medskip

\section{Proof of Theorems \ref{thm2},  \ref{thm3} }\label{Sect6}

\begin{lem}\label{lem:mountS}
	Let $\lambda>0$, $-\infty< \mu <+\infty$. Then \eqref{p} has a mountain pass type solution $u_{\lambda,\mu}\in C^{1,\kappa}(\overline{\Omega})$, $\kappa \in (0,1)$ such that $u_{\lambda,\mu}>0$ in $\Omega$ and	$\Phi_{\lambda,\mu}(u_{\lambda,\mu})>0$. 
\end{lem}
\begin{proof} The functional  $\Phi_{\lambda,\mu}$ satisfies the Palais-Smale condition. Indeed, suppose that $(u_n)\subset  W^{1,p}_0 \setminus 0$ is a Palais-Smale sequence, 
		i.e., 
		$
		\Phi_{\lambda,\mu}(u_n)\to c,~~ D\Phi_{\lambda,\mu}(u_n)\to 0.
		$
		By the Sobolev embedding theorem, we have
	\begin{align*}\label{Pr745}
		c+o(1)\|u_n\|_1=&\frac{\gamma-p}{p\gamma}\|u_n\|_1^p+\lambda\frac{\gamma-q}{q\gamma}\|u_n\|_{L^q}^q-\mu\frac{\gamma-\alpha}{\alpha}\|u_n\|_{L^\alpha}^\alpha
		\geq\\
		&\frac{\gamma-p}{p\gamma}\|u_n\|_1^p-|\mu|\frac{\gamma-\alpha}{\alpha}\|u_n\|_{1}^\alpha,
		~~\mbox{as}~~ n\to \infty.
	\end{align*}
	This  implies that $\|u_n\|_1$ is bounded,  and hence,  after choosing a subsequence if
	necessary,  we have $u_n\rightharpoonup u$ weakly in $W_0^{1,p}$, and $u_n \to u$ strongly
	in $L^r(\Omega)$, $1\leq r<p^*$ to some $u \in W^{1,p}_0$.
	Hence  the convergence $D\Phi_{\lambda,\mu}(u_n)(u_n)\to 0$ as $n\to \infty$ implies
	$$
	\limsup_{n\to \infty} \langle -\Delta_p u_n, u_n-u\rangle =0.
	$$
	Thus by $S^+$ property of the $p$-Laplacian operator (see \cite{DrabekMilota}) it follows that $u_n\to u$ 
	strongly in $W_0^{1,p}$, which means that $\Phi_{\lambda,\mu}$ satisfies the Palais-Smale condition.
	
	The functional $\Phi_{\lambda,\mu}$ possesses a mountain pass type geometry for $\lambda>0$, $-\infty< \mu <+\infty$. Indeed,  for each 
	$\lambda>0$, $-\infty< \mu <+\infty$, there exists $c(\lambda,\mu)>0$ such that $(\lambda/q) s^q-(\mu/\alpha) s^\alpha-(1/\gamma) s^\gamma\geq -c(\lambda,\mu)s^\gamma$,
	$\forall s>0$. Therefore,  by the Sobolev embedding theorem we have 
	\begin{equation}\label{DisZero}
		\Phi_{\lambda,\mu}(u)\geq \frac{1}{p}\|u\|^p_1-c(\lambda,\mu)\|u\|^\gamma_{L^\gamma}\geq (\frac{1}{p}-\tilde{c}(\lambda,\mu)\|u\|^{\gamma-p}_1)\|u\|^p_1,
	\end{equation}
	where $\tilde{c}(\lambda,\mu)>0$ does not depend on $u \in W^{1,p}_0$. We thus can find a sufficiently small $\rho>0$ such that $\Phi_{\lambda,\mu}(u)>\delta$ for some $\delta>0$ provided  $\|u\|_1=\rho$. Evidently, 
	$\Phi_{\lambda,\mu}(su) \to -\infty$ as $t \to +\infty$ for any $u \in  W^{1,p}_0\setminus 0$, and thus, there is $w_1 \in W^{1,p}_0$, $\|w_1\|_1>\rho$ such that $\Phi_{\lambda,\mu}(w_1)<0$. Since $\Phi_{\lambda,\mu}(0)=0$,  $\Phi_{\lambda,\mu}$ possesses a mountain pass type geometry.  It easily seen that the same conclusion holds if we replace the function $f(\mu,\lambda, u):=|u|^{\gamma-2}u+\mu|u|^{\alpha-2}u-\lambda |u|^{q-2}u$ by the truncation function: $f^+(\mu,\lambda, u):=f(\mu,\lambda, u)$ if $u\geq 0$, $f^+(u):=0$ if $u< 0$. Thus,  the mountain pass theorem \cite{AmbrRabin}
	provides us the critical point $u_{\lambda,\mu}$ of $\Phi_{\lambda,\mu}(u_{\lambda,\mu})$ such that 
	$$
	\Phi_{\lambda,\mu}(u_{\lambda,\mu})=\hat{\Phi}_{\lambda,\mu}^{m}:=\inf_{\gamma\in P}\sup_{u\in \gamma} \Phi_{\lambda,\mu}(u)>0
	$$ 
	and $u_{\lambda,\mu}\geq 0$. As in the proof of $(1^o)$, Theorem \ref{thm1},  it follows that  $u_{\lambda,\mu} \in C^{1,\kappa}(\overline{\Omega})$ for  $\kappa \in (0,1)$ and  $u_{\lambda,\mu}>0$ in $\Omega$. 	
\end{proof}

\begin{prop}
	If $\lambda \in (0,\lambda^e)$, $\mu \in (-\infty, \mu^{n,+}_\lambda)$, then $u_{\lambda,\mu}$ is a ground state of \eqref{p}. Moreover, $\Phi_{\lambda,\mu}(u_{\lambda, \mu})\to +\infty$, $\|u_{\lambda, \mu}\|_1 \to +\infty~~\mbox{as}~~\mu \to -\infty$.
\end{prop}
\begin{proof}
	Let $\lambda \in (0,\lambda^e)$, $\mu \in (-\infty, \mu^{n,+}_\lambda)$ and $u_{\lambda,\mu}$ be a  mountain pass solution. Then  $u_{\lambda,\mu} \in \mathcal{N}_{\lambda,\mu}$, and in view of $(2^o)$, Theorem \ref{thm1}, 
	$$
	\Phi_{\lambda,\mu}(u_{\lambda,\mu})=\inf_{g\in P}\max_{s\in [01]} \Phi_{\lambda,\mu}(g(s))\leq \inf_{u \in W^{1,p}_0\setminus 0}\max_{s>0}\Phi_{\lambda,\mu}(su)=\hat{\Phi}^3_{\lambda,\mu}\leq \Phi_{\lambda,\mu}(u_{\lambda,\mu})
	$$ 
	where  $P:=\{g \in C([0,1], W^{1,p}_0):g(0)=0,~ g(1)=w_1\}$ with  $w_1 \in W^{1,p}_0$ such that $\|w_1\|_1>\rho$ and  $\Phi_{\lambda,\mu}(w_1)<0$. Hence $\Phi_{\lambda,\mu}(u_{\lambda,\mu}^3)=\hat{\Phi}^3_{\lambda,\mu}=\Phi_{\lambda,\mu}(u_{\lambda,\mu})$. Thus, for $\lambda \in (0,\lambda^e)$, $\mu \in (-\infty, \mu^{n,+}_\lambda)$, any mountain pass type solution $u_{\lambda,\mu}$ is a ground state of \eqref{p}, i.e.,  $u_{\lambda,\mu} \in \mathcal{M}_{\lambda, \mu}$. Furthermore, by \textit{(iv)}, $(2^o)$, Theorem \ref{thm1}, it follows that $\hat{\Phi}^3_{\lambda,\mu}=\Phi_{\lambda,\mu}(u_{\lambda, \mu})\to +\infty$, $\|u_{\lambda, \mu}\|_1 \to +\infty~~\mbox{as}~~\mu \to -\infty$.
\end{proof}

\medskip

\medskip
\begin{lem}\label{lem:zeroconver} 
Let $\lambda \in (0,\lambda^n)$. Then 
	$\Phi_{\lambda,\mu}(u_{\lambda,\mu})\to 0$ and $\|u_{\lambda, \mu}\|_1 \to 0$  as $\mu \to +\infty$.
\end{lem}
\begin{proof}
	The proof is based on the use of the following auxiliary variational problem
		\begin{equation} \label{GminUnst11}
		\tilde{\Phi}^1_{\lambda,\mu}=\min \{\Phi_{\lambda,\mu}(u):~ u \in \mathcal{R}\mathcal{N}_{\lambda,\mu}^{1}\},
	\end{equation}
	where 
		$$
\mathcal{R}\mathcal{N}_{\lambda,\mu}^1:=\{ u \in \mathcal{N}_{\lambda,\mu}:~\mu^{n,-}_\lambda(u)\leq \mu\},
$$
and $\lambda \in (0,\lambda^n)$, $\mu \in (\mu^{n,-}_\lambda,+\infty)$. 	
	\begin{lem}\label{lem:09_101}
		Let $\lambda \in (0,\lambda^n)$ and  $\mu \in (\mu^{n,-}_\lambda,+\infty)$. There exists a minimizer $\tilde{u}^1_{\lambda,\mu} \in \mathcal{R}\mathcal{N}_{\lambda,\mu}^{1}$ of \eqref{GminUnst11} such that $\tilde{\Phi}^1_{\lambda,\mu}=\Phi_{\lambda,\mu}(\tilde{u}^1_{\lambda,\mu}) >0$.
			\end{lem}
	\begin{proof}
Since $\lambda \in (0,\lambda^n)$, the functional $\mu^{n,-}_\lambda(u)$ is well defined on $W^{1,p}_0\setminus 0$. This   implies that $\mathcal{R}\mathcal{N}_{\lambda,\mu}^1\neq \emptyset$ for $\lambda \in (0,\lambda^n)$, $\mu \in (\mu^{n,-}_\lambda,+\infty)$. By the proof of Lemma \ref{lem:mountS}, there exists $\rho>0$ such that $\inf_{\{u:\|u\|_1=\rho\}}\Phi_{\lambda,\mu}(u)>0$, and therefore  $\tilde{\Phi}^1_{\lambda,\mu}> 0$ for  
$\lambda \in (0,\lambda^n)$,  $\mu \in (\mu^{n,-}_\lambda,+\infty)$.

Let $(u_m)_{m=1}^\infty$ be a minimizing sequence of \eqref{GminUnst11}. 
The coerciveness of $\Phi_{\lambda,\mu}$  on $\mathcal{N}_{\lambda,\mu}$ implies that the sequence $(u_m)$ is bounded in $W^{1,p}_0$ and thus,  up to a
subsequence, 
$$
u_m \to \tilde{u}^1_{\lambda,\mu}~~\mbox{strongly in}~~L^r~\mbox{for}~r \in(1,p^*)~\mbox{and weakly in} ~W^{1,p}_0,
$$
for some $\tilde{u}^1_{\lambda,\mu} \in W^{1,p}_0$. It is easily seen  that if $u_m \to \tilde{u}^1_{\lambda,\mu}$ strongly in $W^{1,p}_0$,  then $\tilde{u}^1_{\lambda,\mu}$ is a non-zero minimizer of \eqref{GminUnst11}.

%Hence  by the weak lower-semicontinuity of  $(\mathcal{R}^n_\lambda)'$  we have
%$$
%0=(\mathcal{R}^n_\lambda)'(t_\lambda^{n,-}(\tilde{u}^1_{\lambda,\mu})\tilde{u}^1_{\lambda,\mu})\leq \liminf_{m\to \infty}(\mathcal{R}^n_\lambda)'(t_\lambda^{n,-}(\tilde{u}^1_{\lambda,\mu}) u_m).
%$$
%This means that  	$t_\lambda^{n,-}(\tilde{u}^1_{\lambda,\mu}) \in (t_\lambda^{n,+}(u_m), t_\lambda^{n,-}(u_m))$, $m=1,\ldots$ and therefore,

To obtain a contradiction, suppose that the convergence $u_m \to \tilde{u}^1_{\lambda,\mu}$ in $W^{1,p}_0$
is not strong. Let us show that  $\tilde{u}^1_{\lambda,\mu}\neq 0$. Observe  $\lim_{m\to +\infty}\|u_m\|_1^p=\beta>0$, since $\tilde{\Phi}^1_{\lambda,\mu}> 0$. Thus, if   $\tilde{u}^1_{\lambda,\mu}=0$, then   $0=\lim_{m\to +\infty}(\Phi_{\lambda,\mu})'(u_m)=(1/p)\beta>0$ is a contradiction.    By the weak lower-semicontinuity of $\mu^{n,-}_\lambda(u)$ we have
$\mu^{n,-}_\lambda(\tilde{u}^1_{\lambda,\mu})\leq \liminf_{m\to \infty}\mu^{n,-}_\lambda(u_m)\leq \mu
$,  and  therefore, there exists $s^1_{\lambda, \mu}(\tilde{u}^1_{\lambda,\mu})>0$ such that
\begin{equation*}\label{RN}
	\mu=\mathcal{R}^n_\lambda(s^1_{\lambda, \mu}(\tilde{u}^1_{\lambda,\mu})\tilde{u}^1_{\lambda,\mu})< \liminf_{m\to \infty}\mathcal{R}^n_\lambda(s^1_{\lambda, \mu}(\tilde{u}^1_{\lambda,\mu})u_m).
\end{equation*}
 Hence   $s^1_{\lambda, \mu}(\tilde{u}^1_{\lambda,\mu})<s^1_{\lambda, \mu}(u_m)$, $m=1,2, \ldots$. In view of that $\Phi'_{\lambda, \mu}(su_m)>0$ for  $s \in (0,s^1_{\lambda, \mu}(u_m))$, this implies
\begin{align*}\label{FGG1}
	\Phi_{\lambda,\mu}(s^1_{\lambda, \mu}(\tilde{u}^1_{\lambda,\mu})&\tilde{u}^1_{\lambda,\mu})<\\
	&\liminf_{m\to \infty}\Phi_{\lambda,\mu}(s^1_{\lambda, \mu}(\tilde{u}^1_{\lambda,\mu})u_m)\leq \liminf_{m\to \infty}\Phi_{\lambda,\mu}(s^1_{\lambda, \mu}(u_m)u_m)=\tilde{\Phi}^{1}_{\lambda,\mu},
\end{align*}
which is a contradiction since $s^1_{\lambda, \mu}(\tilde{u}^1_{\lambda,\mu})\tilde{u}^1_{\lambda,\mu} \in {\mathcal{R}\mathcal{N}}_{\lambda,\mu}^1$.  
\end{proof}
	
	\begin{prop}\label{lem:05_11}
		Let $\lambda \in (0,\lambda^n)$. Then  $\Phi_{\lambda,\mu}(\tilde{u}^1_{\lambda,\mu}) \to 0$ as $\mu \to +\infty$.
			\end{prop}
	\begin{proof}
		Corollary  \ref{cor:contPHI} implies that $\Phi_{\lambda,\mu}(\tilde{u}^1_{\lambda,\mu}) $ is monotone decreasing on $(\mu^{n,-}_\lambda,+\infty)$, and therefore, $\Phi_{\lambda,\mu}(\tilde{u}^1_{\lambda,\mu}) \to \delta$ as $\mu \to +\infty$ for some $\delta\in (0,+\infty)$. 
			Assume by contradiction that $\delta>0$. Then, since $\Phi_{\lambda,\mu}'(\tilde{u}^1_{\lambda,\mu}) =0$, we have
				\begin{equation}\label{limlambda12}
				\frac{\delta}{2}<\Phi(\tilde{u}^1_{\lambda,\mu})=\frac{\gamma-p}{\gamma p}\|\tilde{u}^1_{\lambda,\mu}\|_1^p+\lambda\frac{\gamma-q}{\gamma q}\|\tilde{u}^1_{\lambda,\mu}\|_{L^q}^q-\mu\frac{\gamma-\alpha}{\gamma \alpha}\|\tilde{u}^1_{\lambda,\mu}\|_{L^\alpha}^\alpha<\frac{3\delta}{2},
				\end{equation}
				for sufficiently large $\mu$, whence follows by the embedding $W_0^{1,p}\hookrightarrow L^\alpha(\Omega)$
				\begin{equation}\label{limlambda13}
				\frac{\gamma-p}{\gamma p}\|\tilde{u}^1_{\lambda,\mu}\|_1^p-\mu C\frac{\gamma-\alpha}{\gamma \alpha}\|\tilde{u}^1_{\lambda,\mu}\|_{1}^\alpha<\frac{3\delta}{2},
				\end{equation}
				for some positive constant $C$. Hence $(\tilde{u}^1_{\lambda,\mu})$ is bounded in $W_0^{1,p}$. Consequently, there exists a subsequence $\tilde{u}_{\lambda,\mu_j}^1$ such that $\lim_{j\to +\infty}\mu_j=+\infty$ and $\tilde{u}_{\lambda,\mu_j}^1\rightharpoonup \bar{u}$ weakly in $W^{1,p}_0$ and $\tilde{u}_{\lambda,\mu_j}^1\to \bar{u}$ strongly in $L^r(\Omega),~1\leq r<p^*$ as $j \to \infty$ for some $\bar{u} \in W^{1,p}_0$. Observe that \eqref{limlambda12} implies $\|\tilde{u}^1_{\lambda,\mu}\|_{L^\alpha}^\alpha\to 0$ as $\mu \to+\infty$, which implies $\bar{u}= 0$. Passing to the limit in $\Phi_{\lambda,\mu}'(\tilde{u}_{\lambda,\mu_j}^1)=0$ we obtain $\lim_{j\to +\infty}\|\tilde{u}_{\lambda, \mu_j}^{1}\|_1^p=0$, and consequently,  \eqref{limlambda12} implies that $\mu_j\|\tilde{u}_{\lambda,\mu_j}^1\|_{L^\alpha}^\alpha \to 0$.  Hence  
					$$\frac{\delta}{2}\leq\Phi(\tilde{u}^1_{\lambda,\mu})=\lim_{\mu\to \infty}\Phi(\tilde{u}^1_{\lambda,\mu})=0.$$
			Thus $\delta=0$ and we obtain $\|\tilde{u}^1_{\lambda,\mu}\|_1 \to 0$  and  $\Phi_{\lambda,\mu}(\tilde{u}^1_{\lambda,\mu}) \to 0$ as $\mu \to +\infty$. 
\end{proof}
Let us now conclude the proof of Lemma \ref{lem:zeroconver}. 
 Since $\mu^{n,-}_\lambda(\tilde{u}^1_{\lambda,\mu})\leq \mu$,  the function $\Phi_{\lambda,\mu}(s\tilde{u}^1_{\lambda,\mu})$ has a unique global maximum point $s=s^1_{\lambda, \mu}(\tilde{u}^1_{\lambda,\mu})=1$.  Take a sufficiently large $s_0>1$  such that  $\Phi_{\lambda,\mu}(s_0\tilde{u}^1_{\lambda,\mu})<0$. Then by the above there exists a mountain pass solution $u_{\lambda,\mu}$ such that 
\begin{equation*}
	\Phi_{\lambda,\mu}(u_{\lambda,\mu})=\hat{\Phi}_{\lambda,\mu}^{m}:=\inf_{\gamma\in P}\sup_{u\in \gamma} \Phi_{\lambda,\mu}(u)>0,
\end{equation*} 
where $P:=\{g \in C([0,1];W^{1,p}_0):g(0)=0,~ g(1)=s_0\tilde{u}^1_{\lambda,\mu}\}$. Note that $\tilde{g} \in P$, where $\tilde{g}=s\tilde{u}^1_{\lambda,\mu}$, $s \in [0,s_0]$. Hence 
$$
	\Phi_{\lambda, \mu}(\tilde{u}^1_{\lambda,\mu})=\sup_{s>0} \Phi_{\lambda,\mu}(s\tilde{u}^1_{\lambda,\mu}) \geq \inf_{\gamma\in P}\sup_{u\in \gamma} \Phi_{\lambda,\mu}(u)=\Phi_{\lambda,\mu}(u_{\lambda,\mu}). 
	$$
for any $\lambda \in (0,\lambda^n)$ and  $\mu \in (\mu^{n,-}_\lambda,+\infty)$. Then by Proposition \ref{lem:05_11},  $\Phi_{\lambda,\mu}(u_{\lambda,\mu})\to 0$, and consequently, $\|u_{\lambda, \mu}\|_1 \to 0$ as $\mu \to +\infty$. 
\end{proof}
This concludes the proof of Theorem \ref{thm2}.

\noindent\textit{Proof of Theorem} \ref{thm3}: 

 The existence of three solutions $u_{\lambda,\mu}^i$, $i=1,2,3$, for  $\lambda \in (0,\lambda^e)$ and $\mu\in[\mu_\lambda^{e,-},\mu_\lambda^{n,-})$ follows from  Theorems \ref{thm1}, \ref{thm2}, where we set $u^1_{\lambda,\mu}:=u_{\lambda,\mu}$, for $\lambda \in (0,\lambda^e)$, $\mu\in[\mu_\lambda^{e,-},\mu_\lambda^{n,-})$. They are distinct since $\Phi_{\lambda,\mu}(u^1_{\lambda,\mu})>0$ by Theorem \ref{thm2}, while $\Phi_{\lambda,\mu}(u^2_{\lambda, \mu})<0$, $\Phi_{\lambda,\mu}(u^3_{\lambda, \mu})<0$ and $\Phi''_{\lambda,\mu}(u^2_{\lambda, \mu})>0$, $\Phi''_{\lambda,\mu}(u^3_{\lambda, \mu})<0$  by Theorem \ref{thm1}. By Theorem  \ref{thm1}, $u^2_{\lambda,\mu}$ is linearly stable while $u^3_{\lambda,\mu}$ is linearly unstable solution. 

Consider
$$
K_{\hat{\Phi}_{\lambda,\mu}^{m}}:=\{u \in W_0^{1,p}:~\Phi_{\lambda,\mu}(u)=\hat{\Phi}_{\lambda,\mu}^{m},~D\Phi_{\lambda,\mu}(u)=0\},
$$
where $\Phi_{\lambda,\mu}$ is replaced by the truncation functional as in the proof of Lemma \ref{lem:mountS}. 
Let us show that for $\lambda \in (0,\lambda^e)$ and $\mu\in[\mu_\lambda^{e,-},\mu_\lambda^{n,-})$, $K_{\hat{\Phi}_{\lambda,\mu}^{m}}$ contains a point $u^1_{\lambda,\mu}$ which is a linearly unstable solution. Indeed, by \eqref{DisZero} it is easily seen that $0 \in W_0^{1,p}$ is a local minimizer of   $\Phi_{\lambda,\mu}(u)$. Furthermore, by the proof of $(1^o)$, Theorem \ref{thm1},   $u^2_{\lambda,\mu}$ is also a local minimizer of $\Phi_{\lambda,\mu}(u)$ and $0=\Phi_{\lambda,\mu}(0)>\Phi_{\lambda,\mu}(u^2_{\lambda, \mu})$. Hence, by the result of  Hofer \cite{hofer}, Pucci, Serrin \cite{PuchiSerrin} it follows that the set 
$K_{\hat{\Phi}_{\lambda,\mu}^{m}}$ 
contains a critical point $u^1_{\lambda,\mu}$ which is not local minimum of $\Phi_{\lambda,\mu}(u)$, and therefore it is a linearly unstable solution.

%%%%%%%%%%%%%%%%%%%%%%%%%%%%%%%%%%%%%%%%%%%%%%%%%%%%%%%%%%%%%%%%%%%%%%%%%%%%%%%%%%%%%%%%%5

\appendix 
\section{Appendix}
The statements below are proved using the approach introduced in \cite{ilEuler}.
\begin{prop}\label{PMon} Let $i=1,2,3$, $\lambda,\mu_a, \mu_b \in \mathbb{R}$, $\mu_b>\mu_a$. Assume that $v_{\lambda,\mu_a}^i, v_{\lambda,\mu_b}^i$  are minimizers of \eqref{GminUnst11} for $i=1$, of \eqref{GminUnst0} for $i=2$ and \eqref{GminUnst1} for $i=3$. Then for sufficiently small $|\mu_b-\mu_a|$ there holds
\begin{align}\label{Nerav}
	-\frac{(\mu_b-\mu_a)(s^i_{\lambda, \mu_a}(v_{\lambda,\mu_b}^i))^\alpha}{\alpha}\|v_{\lambda,\mu_b}^i\|_{L^\alpha}^\alpha&<\Phi_{\lambda,\mu_b}(v_{\lambda,\mu_b}^i) 	-  \Phi_{\lambda,\mu_a}(v_{\lambda,\mu_a}^i)<\nonumber\\
	&- \frac{(\mu_b-\mu_a)(s^i_{\lambda, \mu_b}(v_{\lambda,\mu_a}^i))^\alpha}{\alpha}\|v_{\lambda,\mu_a}^i\|_{L^\alpha}^\alpha.
\end{align}
 \end{prop}
\begin{proof} Proofs for $i=1,2,3$ are similar. As an example we prove for  $i=1$. 
		Evidently,  $s^1_{\lambda, \mu_b}(v_{\lambda,\mu_a}^1)v_{\lambda,\mu_a}^1 \in \mathcal{N}_{\lambda,\mu_b}^1$. 
	Moreover, 
	$\mu^{e,-}_\lambda(v_{\lambda,\mu_a}^1)\leq \mu_a<\mu_b$, and therefore, we have  $s^1_{\lambda, \mu_b}(v_{\lambda,\mu_a}^1)v_{\lambda,\mu_a}^1 \in \mathcal{R}\mathcal{N}_{\lambda,\mu_b}^{1}$. Hence $\Phi_{\lambda,\mu_b}(v_{\lambda,\mu_b}^1)\leq\Phi_{\lambda,\mu_b}(s^1_{\lambda, \mu_b}(v_{\lambda,\mu_a}^1)v_{\lambda,\mu_a}^1 )$, and consequently 
			\begin{align}\label{mu12u2}
		\Phi_{\lambda,\mu_b}(v_{\lambda,\mu_b}^1)\leq 
		\Phi_{\lambda,\mu_a}(s^1_{\lambda, \mu_b}(v_{\lambda,\mu_a}^1)v_{\lambda,\mu_a}^1)-\frac{(\mu_b-\mu_a)}{\alpha} \|s^1_{\lambda, \mu_b}(v_{\lambda,\mu_a}^1)v_{\lambda,\mu_a}^1\|_{L^\alpha}^\alpha.
		\end{align}	
Observe, 	if $\mu_b>\mu_a$,  then
	$	0<s_{\lambda,\mu_b}^1(v_{\lambda,\mu_a}^1)<s_{\lambda,\mu_a}^1(v_{\lambda,\mu_a}^1)$. Thus, since $s\mapsto\Phi_{\lambda,\mu_a}(su_{\lambda,\mu_a}^1)$ is increasing in $[0,s^1_{\lambda,\mu_a}(v_{\lambda,\mu_a}^1)]$, we obtain
	$$
	\Phi_{\lambda,\mu_a}(s^1_{\lambda, \mu_b}(v_{\lambda,\mu_a}^1)v_{\lambda,\mu_a}^1)\leq \Phi_{\lambda,\mu_a}(s^1_{\lambda, \mu_a}(v_{\lambda,\mu_a}^1)v_{\lambda,\mu_a}^1)=\Phi_{\lambda,\mu_a}(v_{\lambda,\mu_a}^1),
	$$ 
	and consequently,  the second inequality in \eqref{Nerav}. 	
	The proof of the first inequality in \eqref{Nerav} is handled in much the same way. 
	\end{proof}

	From Lemmas \ref{lem:22_1}, \ref{lem:09_10}, \ref{lem:09_101}, Proposition \ref{PMon} and using the coerciveness of $\Phi_{\lambda,\mu}$  on $\mathcal{N}_{\lambda,\mu}$  it is not hard  to show 	
	\begin{cor}\label{cor:contPHI}
The functions $\mu \mapsto \Phi_{\lambda,\mu}(\tilde{u}_{\lambda,\mu}^1)$ on  $(\mu^{n,-}_\lambda,+\infty)$ for $\lambda \in (0,\lambda^n)$; $\mu \mapsto \Phi_{\lambda,\mu}(u_{\lambda,\mu}^2)$ on  $(\mu^{n,+}_\lambda,\mu^{n,-}_\lambda)$ $\lambda \in (0,\lambda^e)$, $\mu \mapsto \Phi_{\lambda,\mu}(u_{\lambda,\mu}^3)$ on  $(-\infty, +\infty)$  for $\lambda \in (0,\lambda^e)$  are continuous and monotone decreasing.
	\end{cor}

\begin{prop}\label{ContPh}
				Let  $\mu_0, \mu_m \in (\mu^{n,+}_\lambda, \mu^{n,-}_\lambda)$ $(\mu_0, \mu_m \in (-\infty, +\infty))$, $m=1,2 \ldots$ such that $\mu_m \to \mu_0$ as $m \to +\infty$. Then there exist   a subsequence, which we again denote  by $(\mu_m)$,  and a sequence $u^2_{\lambda,\mu_m}\in \mathcal{M}^2_{\lambda, \mu_m}$ $(u^2_{\lambda,\mu_m}\in \mathcal{M}^3_{\lambda, \mu_m})$ such that $ u^2_{\lambda,\mu_m} \to u^2_{\lambda,\mu_0} $ strongly in $W^{1,p}_0$. 
			\end{prop}
\begin{proof} As an example we prove the proposition in the case $i=2$.
By the above it follows that $\Phi_{\lambda,\mu_n}(v_{\lambda,\mu_m}^2) \to \hat{\Phi}^2_{\lambda,\mu_0}$ as $m \to +\infty$, which easily implies that  $(s_{\lambda,\mu_0}^2(v_{\lambda,\mu_m}^2)v_{\lambda,\mu_m}^2)_{m=1}^\infty$ is a minimizing sequence of \eqref{GminUnst0} for $\mu=\mu_0$. Then from the proof of Lemma \ref{lem:22_1} it follows that  up to a
subsequence, $u^2_{\lambda,\mu_m} \to u^2_{\lambda,\mu_0}$ strongly in $W^{1,p}_0$, for some $u^2_{\lambda,\mu_0} \in \mathcal{M}^2_{\lambda, \mu_0}$. 
\end{proof}

%\section*{Acknowledgements}
%
%Yavdat Il'yasov was partially supported by the Russian Science Foundation (No.22-21-00580).


\begin{thebibliography}{0}





\bibitem{Amann} H. Amann,  Fixed point equations and nonlinear eigenvalue problems in ordered Banach spaces. SIAM review, 18(4), (1976) 620-709.

%\bibitem{Amann2}  H. Amann, Existence and multiplicity theorems for semi-linear elliptic boundary value problems. Mathematische Zeitschrift, 150(3),(1976) 281-295.

\bibitem{AmbrRabin} A. Ambrosetti, P. H. Rabinowitz, Dual variational methods in critical point theory and applications, J. Funct. Anal. 14 (4) (1973) 349-381.
 
%\bibitem{Arancibia} A. O. Arancibia,  Certain Aspects of Problems with Non Homogeneous Reactions. In: Smith F.T., Dutta H., Mordeson J.N. (eds) Mathematics Applied to Engineering, Modelling, and Social Issues. Studies in Systems, Decision and Control, vol 200. Springer, Cham. (2019)
%https:doi.org/10.1007/978-3-030-12232-4_3

\bibitem{Arnold} V.I. Arnol'd,  Catastrophe theory, Springer Science \& Business Media, 2003.

\bibitem{MarcCarlIl} M. L. Carvalho, Y. Il'yasov, C. A. Santos, 
\textit{Separating solutions of nonlinear problems using nonlinear generalized Rayleigh quotients}, Topol. Meth. Nonl. Anal. (to appear) (2021), 1--28. 	arXiv:1906.07759 

%\bibitem{Cohen} D.S. Cohen, Multiple stable solutions of nonlinear boundary value problems arising in chemical reactor theory." SIAM Journal on Applied Mathematics 20, no. 1 (1971): 1-13.

\bibitem{Bartnik} R. Bartnik,  L. Simon,  Spacelike hypersurfaces with prescribed boundary values and mean curvature. Communications in Mathematical Physics, 87(1), (1982) 131-152.

\bibitem{Bebernes} J. Bebernes, D. Eberly, Mathematical Problems from Combustion Theory, Springer-Verlag, 1989. 



\bibitem{BobDrHer} V. Bobkov,  P. Dr\'abek, J. Hern\'andez, Existence and multiplicity results for a class of semilinear elliptic equations. Nonlinear Analysis, 200, (2020), 112017.

\bibitem{Brown} K. J. Brown,  M. M. A. Ibrahim, R. Shivaji,   S-shaped bifurcation curves. Nonlinear Analysis: Theory, Methods \& Applications, 5(5), (1981) 475-486. 
%DOI:10.1016/0362-546X(81)90096-1
 
\bibitem{Coelho}  I. Coelho, C. Corsato,  F. Obersnel, P. Omari,   Positive solutions of the Dirichlet problem for the one-dimensional Minkowski-curvature equation. Advanced Nonlinear Studies, 12(3), (2012) 621-638.

%\bibitem{Smoller} J. Smoller, Shock waves and reaction—diffusion equations (Vol. 258). Springer Science \& Business Media. 2012.
\bibitem{CrandRabin} M.G. Crandall, P. H. Rabinowitz, Bifurcation, perturbation of simple eigenvalues,
and linearized stability. Archive for Rational Mechanics and Analysis, 52(2), (1973) 161-180.
%DOI:10.1007/BF00282325 15

\bibitem{cohen}  D.S. Cohen, Multiple stable solutions of nonlinear boundary value problems arising in chemical reactor theory. SIAM Journal on Applied Mathematics, 20(1), (1971) 1-13.

\bibitem{Dancer} E.N. Dancer, On the structure of solutions of an equation in catalysis theory when a parameter is large. Journal of Differential Equations, 37(3),( 1980), 404-437.

\bibitem{DiazHerYA}  J.I.~D\'iaz, J.~Hern\'andez,  Y. Sh.~Il'yasov,
On the exact multiplicity of stable ground states of non-Lipschitz semilinear elliptic equations for some classes of starshaped sets. Advances in Nonlinear Analysis 9, no. 1 (2019) 1046-1065.

\bibitem{DiBened} E. DiBenedetto, $C^{1+\alpha}$ local regularity of weak solutions of degenerate elliptic equations, Nonlinear Anal. 7 (8) (1983) 827–
850. 



\bibitem{DrabekMilota} P. Dr\'abek, J. Milota, Methods of nonlinear analysis, 2nd ed., Birkhauser
Advanced Texts: Basler Lehrbucher. [Birkhauser Advanced Texts: Basel Textbooks], Birkhauser/Springer Basel AG, Basel, 2013.
%MR3025694

\bibitem{gelfand}  I.M. Gel'fand, Some problems in the theory of quasi-linear equations, Uspekhi Mat. Nauk, 14:2(86) (1959), 87–158

\bibitem{Gilmore} R. Gilmore, Catastrophe theory for scientists and engineers. Courier Corporation, 1993.

\bibitem{hofer} H. Hofer, A Geometric Description of the Neighbourhood of a Critical Point Given by the Mountain‐Pass Theorem. Journal of the London Mathematical Society, 2(3) (1985) 566-570.

\bibitem{huang} S.Y. Huang, Classification and evolution of bifurcation curves for the one-dimensional Minkowski-curvature problem and its applications." Journal of Differential Equations 264, no. 9 (2018): 5977-6011.
%https://doi.org/10.1007/s10884-021-09982-4

\bibitem{ilyaReil} Y. Ilyasov,  On extreme values of Nehari manifold
method via nonlinear Rayleigh's quotient, Topol. Methods Nonlinear Anal. 49  (2) (2017) 683-714. 

\bibitem{ilConC} Y. Il'yasov,  On nonlocal existence results for elliptic equations with convex-concave nonlinearities, Nonlinear Anal. 61(1-2)(2005) 211-236.

\bibitem{ilEuler} Y.S. Il'yasov, Euler's Functional for Equations Involving the p-Laplacian as a Function of the Spectral Parameter. Proceedings of the Steklov Institute of Mathematics-Interperiodica Translation, 214 (1996) 175-186.

\bibitem{Kolmog} A. Kolmogoroff, I. Petrovsky, N. Piscounoff, Study of the diffusion equation with
growth of the quantity of matter and its application to a biological problem. (French) Moscow Univ. Bull. Math. 1, (1937) 1-25 

\bibitem{Korman}  P. Korman, Y. Li, On the exactness of an S-shaped bifurcation curve. Proceedings of the American Mathematical Society, 127(4), (1999), 1011-1020.


\bibitem{Lee}  E. Lee, S. Sasi,  R. Shivaji,   S-shaped bifurcation curves in ecosystems. Journal of mathematical analysis and applications, 381(2), (2011), 732-741.

%\bibitem{Kuzin} I.A. Kuzin, S.I. Pohozaev, Entire solutions of semilinear elliptic equations (Vol.33). Birkhäuser. 2012

\bibitem{Lieber}   G.M. Lieberman, Boundary regularity for solutions of degenerate elliptic equations, Nonlinear Anal. 12 (11) (1988) 1203–
1219.

\bibitem{Lubyshev} V. Lubyshev,   Three positive solutions of a nonlinear Dirichlet problem with competing power
nonlinearities, (2015). arXiv:1501.03870


\bibitem{Ludw} Ludwig, D., D. G. Aronson, and H. F. Weinberger. "Spatial patterning of the spruce budworm." Journal of Mathematical Biology 8, no. 3 (1979): 217-258.

%\bibitem{Marsden} J.E. Marsden,  F.J. Tipler, Maximal hypersurfaces and foliations of constant mean curvature in general relativity. Physics Reports, 66(3), (1980), 109-139.

%\bibitem{Mizoguchi} N. Mizoguchi and T. Suzuki, Equations of gas combustion: s-shaped bifurcation and mushrooms,
%J. Differential Equations 134, 183-215 (1997).
 





%
%\bibitem{payne} L. E. Payne,  D. H. Sattinger,  Saddle points and instability of nonlinear hyperbolic equations, Israel J. Math. 22 (3-4) (1975) 273-303.
\bibitem{Poh} S. I. Pokhozhaev, On the method of fibering a solution in nonlinear boundary value problems, Proc. Steklov Inst. Math. 192 (1992) 157-173. 

\bibitem{PuchiSerrin} P. Pucci,  J. Serrin, The structure of the critical set in the mountain pass theorem. Transactions of the American Mathematical Society, 299(1), (1987) 115-132.

\bibitem{Rinzel}  J. Rinzel, W. C. Troy, A one-variable map analysis of bursting in the Belousov-Zhabotinskii reaction. In Oscillations in mathematical biology, 1-23. Springer, Berlin, Heidelberg, 1983.

\bibitem{Tolksdorf} P. Tolksdorf, Regularity for a more general class of quasilinear elliptic equations. J. Diff. Eqs. 51 (1984), 126-50.

\bibitem{Smoller} J. Smoller,  Shockwaves and reaction-diffusion equations. Second edition. Grundlehren
der Mathematischen Wissenschaften, 258. Springer-Verlag, NewYork (1994)

\bibitem{Strogatz} S.H. Strogatz, Nonlinear dynamics and chaos with student solutions manual: With applications to physics, biology, chemistry, and engineering. CRC press, 2018.

%\bibitem{struw} M. Struwe, Variational methods. Vol. 31999. Berlin etc.: Springer, 1990.

\bibitem{trudin} N. S. Trudinger,  On Harnack type inequalities and their application to quasilinear elliptic equations, Comm. Pure and Appl. Math. \textbf{20} (4) (1967),  721--747.

\bibitem{Wang}  S.-H. Wang, On S-shaped bifurcation curves, Nonlinear Analysis TMA 22, (1994), 1475-1485.


\bibitem{zeeman}  E.C. Zeeman,  Catastrophe theory: Selected papers. Bull. Amer. Math. Soc, 84, (1978) 1360-1368.

\bibitem{zeidl} E. Zeidler,   Nonlinear functional analysis and its applications: III: variational methods and optimization. Springer Science \& Business Media, 2013.
\end{thebibliography}
\end{document}